\documentclass[11pt]{article}
\usepackage{arxiv_macros}
\usepackage{algpseudocode}
\usepackage{colortbl}
\definecolor{PineGreen}{RGB}{0,110,51}
\usepackage{xspace}
\usepackage{tikz}
\usetikzlibrary{arrows.meta,positioning,shapes.geometric,fit,calc}

\title{\textbf{Compressed Per-Player Local Methods \\for Multiplayer Federated Learning:\\ Exact Convergence via Increment Sharing}}

\begin{document}

\author{%
  \begin{tabular}{@{}c@{\qquad}c@{\qquad}c@{}}
    Bhavika Balani & TaeHo Yoon & Nicolas Loizou \\[0.3em]
    \texttt{bbalani1@jh.edu} &
    \texttt{tyoon7@jh.edu} &
    \texttt{nloizou@jhu.edu}
  \end{tabular}%
}
\date{
Department of Applied Mathematics \& Statistics\\ 
Mathematical Institute for Data Science (MINDS) \\Johns Hopkins University 
}
\maketitle

\begin{abstract}
Multiplayer Federated Learning (MpFL), introduced in \cite{yoon2025multiplayer}, extends classical federated learning to settings where clients are rational players in a game with individual objectives and strategies, seeking an equilibrium. 
Here, a natural baseline algorithm is \algname{PEARL-SGD}, where each player performs several local SGD steps before exchanging strategies through a central server.
However, the original \algname{PEARL-SGD} formulation does not extend cleanly to compressed communication. Directly compressing the players’ strategies can prevent convergence to the equilibrium, even with deterministic gradients. To address this issue, we reformulate \algname{PEARL-SGD} so that players communicate local strategy increments rather than the strategies themselves. Without compression, the two formulations are equivalent, but with compression, increment sharing preserves exact convergence with deterministic gradients, based on which we design the communication-efficient compressed PEARL-SGD (\algname{C-PEARL-SGD}). Beyond compression, we also consider variance reduction as another axis of improvement, introducing \algname{PEARL-SVRG} and its compressed variant \algname{C-PEARL-SVRG}. Under standard assumptions, we establish convergence guarantees for all variants through a unified parametric inequality. Numerical experiments, including distributed spectrum sensing, demonstrate the practical benefits of our methods over vanilla \algname{PEARL-SGD}.
\end{abstract}

\section{Introduction}

Federated Learning (FL) has emerged as a powerful collaborative learning framework in which multiple clients jointly train machine learning models without sharing their local data~\cite{kairouz2021advances, konevcny2016federated,mcmahan2017communication, li2020federated}. 
In classical FL, clients perform local computations on private data and periodically communicate through a central server to optimize a shared global model, often using methods such as Local SGD~\cite{mcmahan2017communication, Stich2019_local, gorbunov2021local, KhaledMishchenkoRichtarik2020_tighter}.
While it is typically assumed that clients in FL work toward a common objective, it can be restrictive in multi-agent systems where participants have distinct or competing objectives. In such settings, the goal is no longer to find a shared minimizer, but an equilibrium of a multiplayer game. Such problems arise naturally in economic competition, network resource allocation, electricity markets, energy systems, and control~\cite{AhmedAgiza1998_dynamics, PanPavel2007_global, SaadHanPoorBasar2012_gametheoretic, YeHu2017_game, KalyvaPsillakis2024_distributed}.

Multiplayer Federated Learning (MpFL) was recently introduced by \cite{yoon2025multiplayer} to model this setting.
In MpFL, clients are viewed as rational players with individual objectives, seeking an equilibrium while keeping their training data decentralized. A natural baseline algorithm for MpFL is \emph{Per-Player Local SGD} (\algname{PEARL-SGD})~\cite{yoon2025multiplayer}, which extends the idea of local updates to the multiplayer game setting. In each communication round, each player performs several local stochastic gradient steps on its own strategy while keeping the strategies of the other players fixed at the latest synchronization point, after which the server collects the updated strategies, concatenates them, and broadcasts the synchronized joint strategy back to all players.

Local updates reduce the frequency of communication, but synchronization is still a significant overhead because the full joint strategy has to be shared. 
A natural way to further reduce this communication cost is to compress the communicated information, as is commonly done in distributed~\cite{alistarh2017qsgd,lin2018deep,karimireddy2019error,seide20141, stich2018sparsified,condat2022ef} and federated optimization~\cite{reisizadeh2020fedpaq,haddadpour2021federated}. However, \algname{PEARL-SGD} synchronizes the players' updated strategies, and directly compressing the strategy vectors can prevent convergence, even when deterministic gradients are used (\cref{fig:compression-scheme-comparison}).
This can be addressed by letting the players instead communicate the strategy \emph{increments}, i.e., changes in local strategies from the local descent loops, rather than the updated strategies themselves. Mathematically, without compression, the strategy-sharing and increment-sharing formulations of \algname{PEARL-SGD} are equivalent, producing the exact same iterates. However, the distinction is critical when compression is applied; it turns out that increment sharing is a correct reformulation allowing the compressed variant \algname{C-PEARL-SGD} to converge to equilibrium in the deterministic setting.

In this work, we take this new formulation of \algname{PEARL-SGD} as a starting point and design more advanced \algname{PEARL}-type algorithms. We provide a unified framework for convergence analysis that incorporates both compression and standard variance reduction techniques. Consequently, we show that exact convergence is achievable under both compression errors and gradient stochasticity, significantly reducing the communication overhead and improving the learning dynamics of MpFL.

\subsection{Main contributions}
Our main contributions are summarized as follows.

\noindent $\diamond$ \textbf{Exact Convergence via Increment Sharing.} 
We introduce an increment-sharing reformulation of \algname{PEARL-SGD}, and show that it is equivalent to the original formulation without compression (\cref{proposition:pearl-sgd-equivalence}). On the other hand, we show that the increment-sharing variant allows exact convergence to the equilibrium with compression in the deterministic setting (\cref{corollary:c-pearl-sgd-rate}).

\noindent $\diamond$  \textbf{Advanced MpFL methods.}
Building on increment sharing, we design new per-player local (\algname{PEARL}) methods for MpFL. We develop a generalized \algname{PEARL-SGD} (\cref{alg:unified-pearl}) algorithm and a unified convergence framework that allows different local gradient estimators and compression maps. Consequently, we analyze \algname{PEARL-SGD}, \algname{PEARL-SVRG} and their compressed variants \algname{C-PEARL-SGD} and \algname{C-PEARL-SVRG} in a unified manner (all as corollaries of the main result \cref{thm:main_theorem}), the latter three being new methods.

\noindent $\diamond$  \textbf{Numerical evaluation.}
We complement our theory with experiments showing practical effects of local updates, compression, and variance reduction in MpFL. 

\section{Multiplayer Federated Learning and Generalized \texorpdfstring{\\ \algname{PEARL-SGD}}{PEARL-SGD}}
\label{section:MpFL-concept}

In this section, we review the MpFL setting and the \algname{PEARL-SGD} algorithm from \cite{yoon2025multiplayer}, and prove its equivalence to our increment-sharing variant. We then introduce a generalized \algname{PEARL-SGD} framework that incorporates variance reduction and compression.\smallskip

\noindent
\textbf{Multiplayer games and equilibrium.}
We begin by introducing the concept of a multiplayer game and the corresponding notion of equilibrium. Throughout this work, players are indexed by $i\in [n]$, and their respective strategies are denoted by $x^i \in \mathbb{R}^{d_i}$.
Let $\mathbf{x} = (x^1, \dots , x^n) \in \mathbb{R}^D$,  with $D=\sum_{i=1}^n d_i$, be the joint strategy vector of all players. Each player $i\in[n]$ tries to minimize its own objective function $f_i(x^i; x^{-i}) = f_i(x^1, \dots, x^n)\colon \mathbb{R}^{D} \to \mathbb{R}$ with respect to its strategy $x^i \in \mathbb{R}^{d_i}$, given all the other players' strategies $x^{-i} = (x^1, \dots, x^{i-1}, x^{i+1}, \dots, x^n) \in \mathbb{R}^{D - d_i}$.
In our setting, the goal is to find an \textit{equilibrium} $\mathbf{x}_\star = (x_\star^1, \dots,x_\star^n ) \in \mathbb{R}^{D}$ satisfying
\begin{align}
\label{eqn:equilibrium}
\begin{split}
    f_i(x_\star^i; x_\star^{-i}) \le f_i (x^i; x_\star^{-i}), \quad
    \forall x^i \in \reals^{d_i}, \quad \forall i \in [n] .
\end{split}
\end{align}
We consider stochastic objectives of the form $f_i(x^i;x^{-i}) = \mathbb E_{\xi^i\sim\mathcal D_i}\left[f_{i,\xi^i}(x^i;x^{-i})\right].$ Here $\mathcal D_i$ denotes the data distribution of player $i$, and $f_{i,\xi^i}$ is its loss for a data point $\xi^i$ sampled from $\mathcal D_i$. For the variance-reduced methods, we additionally assume that each $f_i(x^i; x^{-i}) = \frac{1}{M_i}\sum_{m=1}^{M_i} f_{i,m}(x^i; x^{-i})$ has a finite-sum structure and that player $i$ has access to stochastic gradients of $f_i$ via uniform sampling from the data indices $1,\dots,M_i$. 

The algorithms considered in this paper may involve two sources of randomness: the mini-batch $\xi_k^i$ used to construct the local gradient estimator $g_k^i$ and, when compression is used, the compressor $\mathcal Q^i$ applied to the local increment $u_p^i$. We assume that the mini-batches are sampled independently of the previous history, and that the compressor randomness is independent of both the previous history and the mini-batch sampling. We use $\mathbb E_{\xi_k^i}$ for expectation over the mini-batch sampling, $\mathbb E_{\mathcal Q^i}$ or $\mathbb E_{\mathcal Q}$ for expectation over the compressor randomness conditioned on the vector being compressed, and $\mathbb E[\cdot]$ for the total expectation over all randomness in the algorithm. Throughout the paper, for simplicity of notation, $\nabla f_i(x^i;x^{-i})$ denotes the partial gradient $\nabla_{x^i}f_i(x^i;x^{-i})$, where we differentiate only with respect to $x^i$ and not $x^{-i}$. We use the same convention for component and mini-batch gradients.\medskip

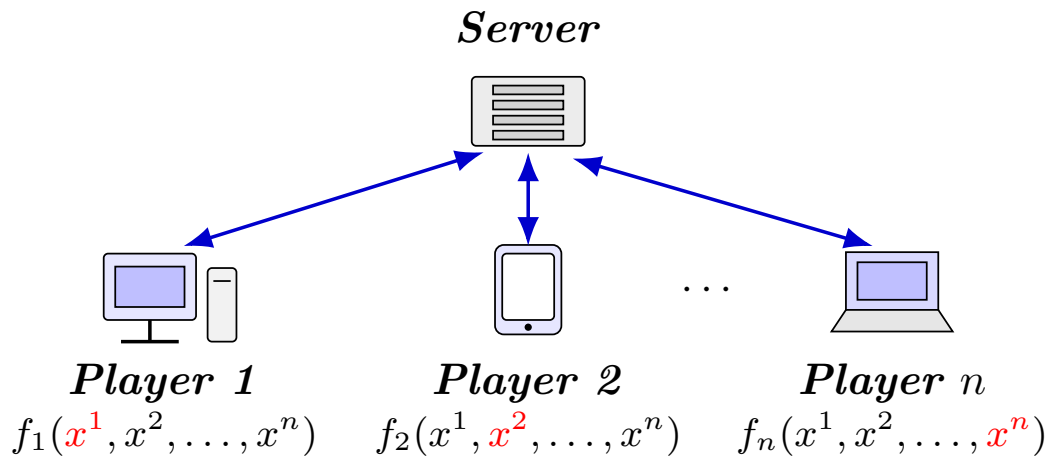
\begin{wrapfigure}{r}{0.43\columnwidth}
\centering
\vspace{-10pt}
\resizebox{0.94\linewidth}{!}{%
\begin{tikzpicture}[
    font=\small,
    >=Latex,
    comm/.style={<->, thick, blue!80!black},
    playerblock/.style={align=center,inner sep=0.3pt},
]

\node[align=center] at (0,2.95) {\textbf{\textit{Server}}};
\draw[rounded corners=1pt, fill=gray!15, draw=black] (-0.45,2.00) rectangle (0.45,2.55);
\foreach \y in {2.05,2.17,2.29,2.41} {
\draw[fill=gray!35, draw=black] (-0.28,\y) rectangle (0.28,\y+0.07);
}

\coordinate (P1) at (-2.90,0.75);
\coordinate (P2) at (0,0.75);
\coordinate (Pn) at (2.90,0.75);

\begin{scope}[shift={(P1)},scale=1.15]
\draw[rounded corners=1pt, fill=blue!10, draw=black] (-0.42,-0.12) rectangle (0.22,0.34);
\draw[fill=blue!25, draw=black] (-0.34,-0.03) rectangle (0.14,0.25);
\draw[thick] (-0.10,-0.12) -- (-0.10,-0.27);
\draw[thick] (-0.30,-0.27) -- (0.10,-0.27);
\draw[rounded corners=1pt, fill=gray!10, draw=black] (0.30,-0.27) rectangle (0.50,0.24);
\draw (0.34,0.15) -- (0.46,0.15);
\end{scope}

\begin{scope}[shift={(P2)},scale=1.15]
\draw[rounded corners=2pt, fill=blue!10, draw=black] (-0.23,-0.22) rectangle (0.23,0.40);
\draw[rounded corners=1pt, fill=white, draw=black] (-0.18,-0.12) rectangle (0.18,0.34);
\fill (0,-0.17) circle (0.025);
\end{scope}

\begin{scope}[shift={(Pn)},scale=1.15]
\draw[fill=blue!15, draw=black] (-0.32,-0.05) rectangle (0.32,0.35);
\draw[fill=blue!25, draw=black] (-0.24,0.02) rectangle (0.24,0.28);
\draw[fill=gray!20, draw=black] (-0.42,-0.20) -- (0.42,-0.20) -- (0.32,-0.05) -- (-0.32,-0.05) -- cycle;
\end{scope}

\draw[comm] (-0.35,1.95) -- ($(P1)+(0.15,0.45)$);
\draw[comm] (0,1.95) -- ($(P2)+(0,0.45)$);
\draw[comm] (0.35,1.90) -- ($(Pn)+(-0.15,0.45)$);
\node at (1.45,0.85) {$\cdots$};

\node[playerblock] at ($(P1)+(0,-0.85)$) {\textbf{\textit{Player 1}}\\[-1.2pt] {\footnotesize $f_1({\color{red}x^1},x^2,\ldots,x^n)$}};
\node[playerblock] at ($(P2)+(0,-0.85)$) {\textbf{\textit{Player 2}}\\[-1.2pt] {\footnotesize $f_2(x^1,{\color{red}x^2},\ldots,x^n)$}};
\node[playerblock] at ($(Pn)+(0,-0.85)$) {\textbf{\textit{Player $n$}}\\[-1.2pt] {\footnotesize $f_n(x^1,x^2,\ldots,{\color{red}x^n})$}};

\end{tikzpicture}%
}

\caption{{\footnotesize Illustration of the MpFL setting. Each player $i$ controls its own strategy $x^i$ and has an individual objective $f_i(x^1,\ldots,x^n)$. The players exchange information through a central server with the goal of reaching an equilibrium $\mathbf{x}_\star\in\mathbb R^D$.}}
\label{fig:mpfl}
\vspace{-10pt}
\end{wrapfigure}
\vspace{.1cm}

\noindent
\textbf{Multiplayer federated learning (MpFL).} 
The MpFL framework focuses on solving an $n$-player game (find an equilibrium point~\eqref{eqn:equilibrium}) when the players communicate via a central server to reach an equilibrium. Players exchange their strategy vectors $x^i$ only when \textit{communication/synchronization} through the central server occurs. Between synchronizations, each player adjusts their strategies based on the recently received information of the other players' strategies through local optimization loops. These parallel local computations are called a single \textit{communication round}. 

\vspace{.2cm}
\noindent
\textbf{\algname{PEARL-SGD} algorithm.}
In order to find $\vx_\star$ satisfying \eqref{eqn:equilibrium} in a communication-efficient manner, \cite{yoon2025multiplayer} proposed
\algname{PEARL-SGD}, where each player runs a fixed number of iterations $\tau$ within each communication round. At the beginning of the $p$-th round, a synchronization step occurs in which the server collects all players' current strategies $x_{\tau p}^i$ ($i\in [n]$) and then concatenates and distributes the joint strategy vector $\vx_{\tau p}$ to all players. Each player $i \in [n]$ then locally runs $\tau$ SGD steps to update their strategy $x^i$ with respect to
their local objective function $f_i$, while keeping the strategies $x^{-i}$ of other players obtained from the previous synchronization step fixed.
After these local updates, the algorithm proceeds to the next round. 

\begin{figure}[t]
\centering
\includegraphics[width=0.45\textwidth]{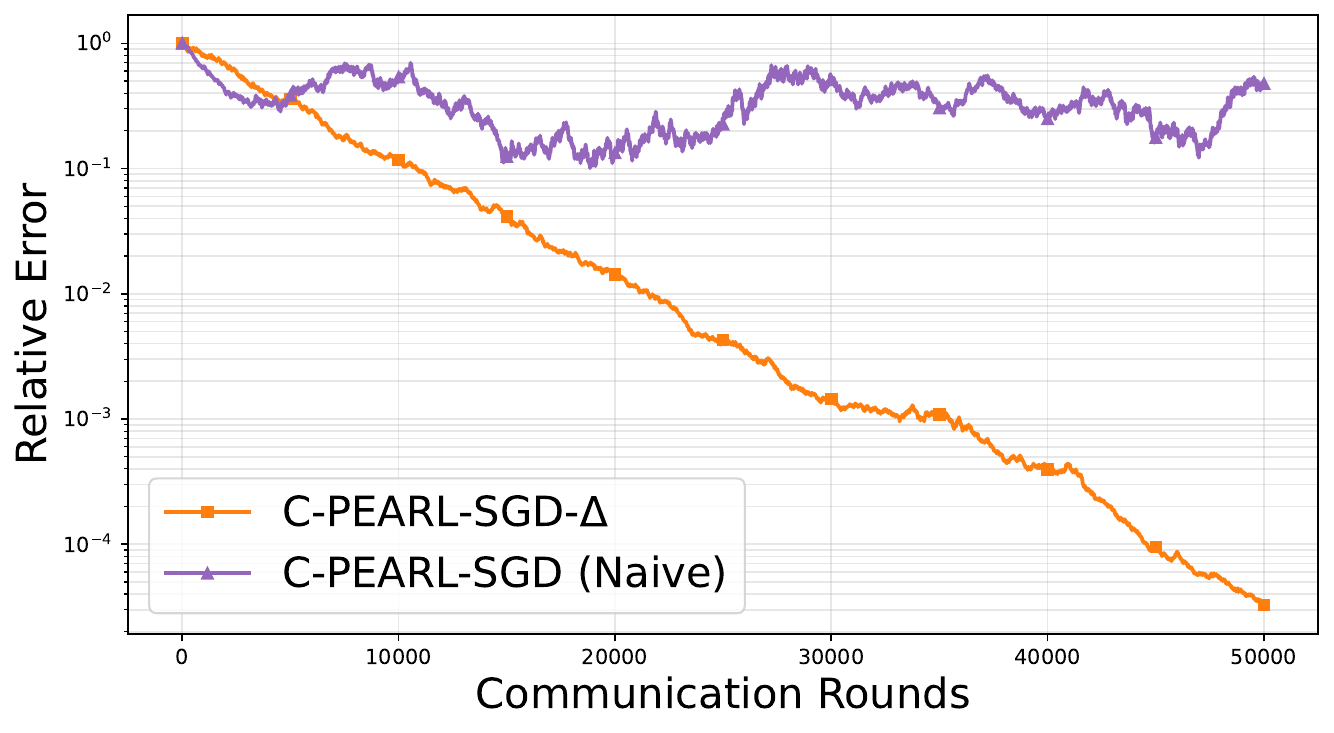}%
\caption{Comparison of convergence for versions of \algname{C-PEARL-SGD} compressing $x_{\tau p}^i$ (Naive) vs.\ $x_{\tau(p+1)}^i - x_{\tau p}^i$ ($\Delta$) in the $n$-player quadratic game setting of \cref{section:experiments} with deterministic gradients.}
\label{fig:compression-scheme-comparison} 
\end{figure}

One limitation of \algname{PEARL-SGD} is that the full joint strategy vector $\vx$ must be transmitted to every player at each synchronization. Unlike classical FL, where local vectors from clients are averaged at the server and only the average is distributed back to the clients, direct communication of $\vx$ to all players in MpFL incurs significant communication overhead. This motivates us to equip each player $i \in [n]$ with compressor maps $\cQ^i$ and communicate the vectors compressed by $\cQ^i$ rather than the full information. However, when we directly apply compression to the strategy vectors $x_{\tau p}^i$ shared in the original \algname{PEARL-SGD} formulation, it may not converge even for well-conditioned deterministic problems, as illustrated by the ``Naive'' version of \textit{Compressed PEARL-SGD} (\algname{C-PEARL-SGD}) in \cref{fig:compression-scheme-comparison}.

Instead, we use an increment-sharing reformulation of \algname{PEARL-SGD}. In this form, after the local training round \(p\), player \(i\) communicates the local increment $u_p^i = \widetilde x_{\tau(p+1)}^i-\widetilde x_{\tau p}^i,$ rather than the updated strategy $\widetilde x_{\tau(p+1)}^i$. Motivated by the empirical behavior of \algname{C-PEARL-SGD} ($\Delta$) in \cref{fig:compression-scheme-comparison}, we introduce the Generalized \algname{PEARL-SGD} (Algorithm \ref{alg:unified-pearl}) update rule based on this increment-sharing mechanism and provide its convergence analysis.

\subsection{Equivalence between strategy and increment sharing}
\label{nasodal}

Let us consider the special case of the generalized \algname{PEARL-SGD}(Algorithm \ref{alg:unified-pearl}) where $\cQ^i=\mathrm{Id}$ (no compression) and $g_k^i$ is the mini-batch gradient 
\begin{align}
    \label{eq:mini_batch_gradient}
    g^i_k = \nabla f_{i,\xi_k^i}(\widetilde x_k^i;x_{\tau p}^{-i}) := \frac{1}{\left| \xi_k^i \right|} \sum_{m \in \xi_k^i} \nabla f_{i,m} (\widetilde x_k^i;x_{\tau p}^{-i}).
\end{align}
Here $\xi_k^i$ denotes a mini-batch of data points sampled from $\mathcal D_i$. In the finite-sum setting, its elements are sampled component indices. The resulting no-compression increment-sharing algorithm, which we refer to as \algname{PEARL-SGD-$\Delta$}, differs from the original form of \algname{PEARL-SGD} from~\cite{yoon2025multiplayer} only in what is communicated. In \algname{PEARL-SGD}, the vectors sent from player $i$ to the master server and from the server to each player are respectively $x_{\tau(p+1)}^i$ and $\vx_{\tau(p+1)}$, rather than the update vectors $q_p^i = u_p^i = \widetilde x_{\tau(p+1)}^i - \widetilde x_{\tau p}^i$ and $\vq_p$ used in \algname{PEARL-SGD-$\Delta$}. The following proposition states that this difference does not change the iterates generated by the algorithm.

\begin{proposition}
\label{proposition:pearl-sgd-equivalence}
\algname{PEARL-SGD-$\Delta$}, obtained by letting $\cQ^i = \mathrm{Id}$ and $g_k^i = \nabla f_{i,\xi_k^i} (\widetilde x_k^i;x_{\tau p}^{-i})$ in Algorithm \ref{alg:unified-pearl}, is equivalent to \algname{PEARL-SGD} in the sense that, with the same initialization,
step-sizes, synchronization interval, and mini-batch realizations, they generate the same sequence of local strategy vectors $\widetilde  x_k^i$ for all $k=0,\dots,\tau R$ and $i=1,\dots,n$.
\end{proposition}
\begin{proof}
We use induction on $k$ to show that $x_k^i$ from the two algorithms are the same for all $i=1,\dots,n$.
For $k=0$, there is nothing to prove.
Now let $0 \le p \le R-1$, and assume the induction hypothesis holds for all $\ell=0,\dots,\tau p$ and that $\widetilde{x}_{\tau p}^i = x_{\tau p}^i$.
Then, for \algname{PEARL-SGD-$\Delta$} we have $\widetilde x_{k+1}^i = \widetilde x_k^i - \gamma_k g_k^i$ for all $k=\tau p,\dots,\tau (p+1)-1$, and finally,
\[
    x_{\tau (p+1)}^i = x_{\tau p}^i + q_p^i = \widetilde x_{\tau p}^i + \left( \widetilde x_{\tau (p+1)}^i - \widetilde x_{\tau p}^i \right) = \widetilde x_{\tau (p+1)}^i
\]
because $q_p^i = \cQ^i (u_p^i) = \widetilde x_{\tau (p+1)}^i - \widetilde x_{\tau p}^i$.
Thus for each $k=\tau p + 1, \dots , \tau (p+1)$, we have $\widetilde x_k^i = \widetilde x_{\tau p}^i - \sum_{\ell=\tau p}^{k-1} \gamma_\ell g_\ell^i$ and thus they coincide with \algname{PEARL-SGD} strategy vectors.
\end{proof}
For the remainder of the paper, we refer to \algname{PEARL-SGD-$\Delta$} simply as \algname{PEARL-SGD}.

\subsection{Exact convergence via increment sharing}
The comparison in~\cref{fig:compression-scheme-comparison} raises a natural question: why does compressing local increments work while directly compressing strategies might not? In this section, we provide an intuitive explanation of this difference. 
We start by stating a standard condition for compressors in the literature~\cite{alistarh2017qsgd,gorbunov2021local,horvath2023stochastic}.
\begin{assumption}
    \label{assumption:unbiased-compression}
A randomized operator $\mathcal{Q}\colon \reals^D \to \reals^D$ is an unbiased compressor with parameter $\omega \ge 0$ if, for every $v\in\reals^D$, $\mathbb{E}\left[\mathcal{Q}(v)\,\middle|\,v\right] = v,$ and \[\mathbb{E}\left[\sqnorm{\mathcal{Q}(v)-v}\,\middle|\,v\right]\le \omega \sqnorm{v}\] where $\mathcal{Q}(v) =\bigl(\mathcal{Q}^1(v^1),\dots,\mathcal{Q}^n(v^n)\bigr)$ applies independent compression across players.
\end{assumption}

Here $\omega$ controls the compression noise and $\omega=0$ corresponds to the case of full communication, where $\cQ = \mathrm{Id}$. 
Note that the expected (squared) norm of the compression error at $v$ tends to zero, as $v$ itself vanishes. 
As the players approach the equilibrium, assuming that deterministic gradients are used, the local updates will produce nearly zero increments in strategies and compression errors become commensurately small. This is not the case if strategy vectors are directly compressed.

Analogously, in distributed optimization, quantized SGD \cite{alistarh2017qsgd} may not converge to the global minimizer $x_*$ even in the deterministic setting, as local gradients may not vanish at $x_*$.
This is addressed by DIANA~\cite{mishchenko2025distributed}, which incorporates an additional memory vector to make the compressed quantity approach zero as the algorithm trajectory approaches $x_*$.
The increment-sharing reformulation of \algname{PEARL-SGD} is analogous to DIANA in the sense that we select the right quantity to compress; however, it does not require an additional memory sequence as in DIANA.
We provide a more detailed comparison in Appendix~\ref{appendix:details_exact_convergence_increment_sharing}. 

\subsection{Assumptions for convergence analysis}

For our convergence analysis, we use the same assumptions on the multiplayer game and the joint gradient operator as in the original analysis of vanilla \algname{PEARL-SGD}~\cite{yoon2025multiplayer}, which are standard in stochastic games and variational inequalities~\cite{LoizouBerardGidelMitliagkasLacoste-Julien2021_stochastic, beznosikov2020distributed, choudhury2024single, gorbunov2022stochastic, GorbunovLoizouGidel2022_extragradient}. For the variance-reduced methods, \algname{PEARL-SVRG} and \algname{C-PEARL-SVRG}, we additionally require the finite-sum structure introduced above and component-wise smoothness (Assumption~\ref{assumption:finite-sum-smoothness}), as is standard for SVRG-type methods.

\paragraph{Assumptions on the multiplayer game}
For every player $i\in[n]$, we assume that $f_i(\cdot;x^{-i})$ is convex and smooth with respect to player $i$'s strategy, while the other players' strategies are fixed.

\begin{assumption}[\textbf{\emph{Convexity (CVX)}}]
\label{assumption:convexity}
For every $i\in[n]$, $x^{-i}\in\mathbb R^{D-d_i}$, and $x^i,y^i\in\mathbb R^{d_i}$,
\vspace{-0.5em}
\begin{align*}
    f_i(y^i;x^{-i})\geq f_i(x^i;x^{-i})+\left\langle\nabla f_i(x^i;x^{-i}),y^i-x^i\right\rangle.
\end{align*}
\end{assumption}

\begin{assumption}[\textbf{\emph{Smoothness (SM)}}]
\label{assumption:smoothness}
For every $i\in[n]$, $x^{-i}\in\mathbb R^{D-d_i}$, and $x^i,y^i\in\mathbb R^{d_i}$,
\vspace{-0.5em}
\begin{align*}
    \left\|\nabla f_i(x^i;x^{-i})-\nabla f_i(y^i;x^{-i})\right\| \leq L_i\left\|x^i-y^i\right\|.
\end{align*}
\end{assumption}

For the non-variance-reduced methods, we additionally assume bounded variance of the stochastic gradients.

\begin{assumption}[\textbf{\emph{Unbiased stochastic gradients and bounded variance (BV)}}]
\label{assumption:noise}
For every player \(i\in[n]\), let \(\xi^i\sim\mathcal D_i\) denote the local sampling randomness. For any \(x^i\in\mathbb R^{d_i}\) and \(x^{-i}\in\mathbb R^{D-d_i}\), the stochastic gradient satisfies $ \mathbb E_{\xi^i}\left[\nabla f_{i,\xi^i}(x^i;x^{-i})\right]= \nabla f_i(x^i;x^{-i}),$
and there exists \(\sigma_i\ge0\) such that $\mathbb E_{\xi^i}\left[\left\|\nabla f_{i,\xi^i}(x^i;x^{-i})-\nabla f_i(x^i;x^{-i})\right\|^2
\right]\le\sigma_i^2.$
\end{assumption}

\paragraph{Assumptions on the joint gradient operator}
We now state the assumptions on the joint gradient operator $F\colon\mathbb R^D\to\mathbb R^D$, defined by $F(\vx)=\left(\nabla f_1(x^1;x^{-1}),\dots,\nabla f_n(x^n;x^{-n})\right)$. Our analysis is carried out around an equilibrium $\vx_\star=(x_\star^1,\dots,x_\star^n)$ satisfying $F(\vx_\star)=0$. We assume that \(F\) satisfies quasi-strong monotonicity and star-cocoercivity. These conditions are weaker than imposing strong monotonicity and cocoercivity globally over all pairs of points~\cite{LoizouBerardGidelMitliagkasLacoste-Julien2021_stochastic, beznosikov2023stochastic}. As a result, quasi-strong monotonicity can hold for some non-monotone games, while star-cocoercivity can hold beyond the globally Lipschitz case.

\begin{assumption}[\textbf{\emph{Quasi-strong monotonicity (QSM)}}]
\label{assumption:quasi-strong-monotonicity}
There exists $\mu>0$ such that, for every $\vx\in\mathbb R^D$
\[\left\langle F(\vx),\vx-\vx_\star\right\rangle \geq\mu\left\|\vx-\vx_\star\right\|^2.\]
\end{assumption}

\begin{assumption}[\textbf{\emph{Star-cocoercivity (SCO)}}]
\label{assumption:star-cocoercivity}
There exists $\ell>0$ such that, for every $\vx\in\mathbb R^D$,
\[\left\langle F(\vx),\vx-\vx_\star\right\rangle
\geq\frac{1}{\ell}\left\|F(\vx)\right\|^2.\]
\end{assumption}
Together, \hyperref[assumption:quasi-strong-monotonicity]{\textbf{\emph{(QSM)}}} and \hyperref[assumption:star-cocoercivity]{\textbf{\emph{(SCO)}}} imply $\mu\|\vx-\vx_\star\|\leq\|F(\vx)\|\leq\ell\|\vx-\vx_\star\|$ for all $\vx\in\mathbb R^D$. Hence $\mu\leq\ell$, and $\kappa=\ell/\mu\geq1$ is the \textit{condition number} of the game. 

\section{Advanced \algname{PEARL} Variants for MpFL}
\label{section:unified_framework_and_assumptions}
Algorithm~\ref{alg:unified-pearl} provides two directions of flexibility: the local gradient estimator \(g_k^i\) and the compression map \(\mathcal Q^i\). Different choices of \(g_k^i\) and \(\mathcal Q^i\) give the \algname{PEARL} variants considered in this work. In this section, we introduce these variants as special cases of generalized \algname{PEARL-SGD} (Algorithm~\ref{alg:unified-pearl}) and discuss their benefits.

\begin{algorithm}[t]
\caption{Generalized \algname{PEARL-SGD}}
\label{alg:unified-pearl}
\algrenewcommand\alglinenumber[1]{}
\begin{algorithmic}
\State \textbf{Input:} Step-size \(\gamma_k >0\), synchronization interval \(\tau\geq1\), communication rounds \(R\geq1\), 
compressor maps \(\{\mathcal Q^i\}_{i=1}^n\), initial local strategies \(\{x_0^i\}_{i=1}^n\)
\vspace{.3cm}
\State Master collects \(x_0^i\) from players \(i\in[n]\)
\State Master broadcasts \(\mathbf{x}_{-\tau}=\mathbf{x}_0\) back to all players
\State \textbf{for} \(p=0,\ldots,R-1\) \textbf{do}
\State \hspace{0.5cm} \textbf{for} \(i=1,\ldots,n\) in parallel \textbf{do}
\State \hspace{1.0cm} Set \(\mathbf{x}_{\tau p}\gets \mathbf{x}_{\tau(p-1)}+\mathbf{q}_{p-1}\), where \(\mathbf{q}_{-1}=0\)
\State \hspace{1.0cm} \(\widetilde{x}_{\tau p}^i\gets x_{\tau p}^i\)
\State \hspace{1.0cm} \textbf{for} \(k=\tau p,\ldots,\tau(p+1)-1\) \textbf{do}
\State \hspace{1.5cm} Sample a mini-batch \(\xi_k^i\sim\mathcal D_i\) 
\State \hspace{1.5cm} Compute an estimator \(g_k^i\)  satisfying: 
 \State \hspace{1.5cm} \(\mathbb E_{\xi_k^i}\left[g_k^i\mid\widetilde x_k^i,x_{\tau p}^{-i}\right]
    =\nabla f_i(\widetilde x_k^i;x_{\tau p}^{-i})\)
\State \hspace{1.5cm} \(\widetilde{x}_{k+1}^i\gets \widetilde{x}_k^i-\gamma_k g_k^i\)
\State \hspace{1.0cm} \textbf{end for}
\State \hspace{1.0cm} \(u_p^i\gets \widetilde{x}_{\tau(p+1)}^i-\widetilde{x}_{\tau p}^i\)
\State \hspace{1.0cm} \(q_p^i\gets \mathcal Q^i(u_p^i)\)
\State \hspace{1.0cm} Player \(i\) sends \(q_p^i\) to master
\State \hspace{0.5cm} \textbf{end for}
\State \hspace{0.5cm} Master collects \(q_p^i\) from players \(i\in[n]\)
\State \hspace{0.5cm} Master broadcasts \(\mathbf q_p=(q_p^1,\ldots,q_p^n)\) to players
\State \textbf{end for}
\State \textbf{Output:} \(\vx_{\tau (R-1)} + \vq_{R-1}  =\mathbf{x}_{\tau R} \in\mathbb R^D\)
\end{algorithmic}
\end{algorithm}

\subsection{Variance Reduction via \algname{PEARL-SVRG}}
\algname{PEARL-SVRG} is a variance-reduced \cite{gorbunov2020unified,palaniappan2016stochastic,yang2020global,ChavdarovaGidelFleuretLacoste-Julien2019_reducing} MpFL method obtained as a special case of Algorithm~\ref{alg:unified-pearl} with $\mathcal Q^i = \mathrm{Id}$ and the variance-reduced gradient estimator
\begin{equation}
\label{eq:pearl-svrg-gradient-estimator}
    g_k^i= \nabla f_{i,\xi_k^i}(\widetilde{x}_k^i;x_{\tau p}^{-i})- \nabla f_{i,\xi_k^i}(x_{\tau p}^i;x_{\tau p}^{-i}) +\nabla f_i(x_{\tau p}^i;x_{\tau p}^{-i}).
\end{equation} 
At the beginning of communication round \(p\), player \(i\) computes the full gradient \(\nabla f_i(x_{\tau p}^i;x_{\tau p}^{-i})\) once at the synchronized strategy \(\mathbf{x}_{\tau p}\) and uses it throughout the \(\tau\) local updates. At each local update, the player computes \(g_k^i\) (using \eqref{eq:pearl-svrg-gradient-estimator}) and updates its local strategy as in Algorithm~\ref{alg:unified-pearl}. Since \(\mathcal Q^i=\mathrm{Id}\), the full strategy increment is communicated at synchronization. \algname{PEARL-SVRG} converges exactly to the equilibrium in the stochastic finite-sum setting, at the cost of computing one full gradient per communication round and two mini-batch gradients per local step. For its convergence analysis, we additionally assume smoothness of each component function $f_{i,m}$.

\begin{assumption}
\label{assumption:finite-sum-smoothness}
For all \(i\in[n]\), \(m\in[M_i]\) and \(x^{-i}\in\reals^{D-d_i}\), the function $f_{i,m}(\cdot;x^{-i})\colon \reals^{d_i}\to\reals$
is \(L_i\)-smooth. That is, for all \(u,v\in\reals^{d_i}\),
\begin{align*}
    \left\|\nabla f_{i,m}(u;x^{-i}) - \nabla f_{i,m}(v;x^{-i})\right\|\leq L_i\|u-v\|.
\end{align*}
\end{assumption}

\begin{table*}[t]
    \centering
    \caption{Special cases recovered by Algorithm~\ref{alg:unified-pearl}. \(\mathrm{Id}\) is the identity map, and \(\mathcal Q^i\) is a general unbiased compressor.}
    \label{tab:unified-pearl-special-cases}
    \scriptsize
    \setlength{\tabcolsep}{2.0pt}
    \renewcommand{\arraystretch}{0.5}
    \resizebox{\textwidth}{!}{%
    \begin{tabular}{c c c}
        \toprule
        \textbf{Method} & \textbf{Unbiased estimator \(g_k^i\)} &  \textbf{Compression Operator}\\
        \midrule
        \algname{PEARL-SGD}
        &
        \(\nabla f_{i,\xi_k^i}
        (\widetilde{x}_k^i;x_{\tau p}^{-i} )\)
        &
        \(\mathrm{Id}\)
        \\
        \algname{PEARL-SVRG}
        &
        \(\begin{aligned}
       \nabla f_{i,\xi_k^i}(\widetilde{x}_k^i;x_{\tau p}^{-i}) - \nabla f_{i,\xi_k^i}(
        x_{\tau p}^i;x_{\tau p}^{-i})+ \nabla f_i(x_{\tau p}^i;x_{\tau p}^{-i})
        \end{aligned}\)
        &
        \(\mathrm{Id}\)
        \\
        \algname{C-PEARL-SGD}
        &
        \(\nabla f_{i,\xi_k^i}
        (\widetilde{x}_k^i;x_{\tau p}^{-i} )\)
        &
        \(\mathcal Q^i\)
        \\
        \algname{C-PEARL-SVRG}
        &
         \(\begin{aligned}
       \nabla f_{i,\xi_k^i}(\widetilde{x}_k^i;x_{\tau p}^{-i}) - \nabla f_{i,\xi_k^i}(
        x_{\tau p}^i;x_{\tau p}^{-i})+ \nabla f_i(x_{\tau p}^i;x_{\tau p}^{-i})
        \end{aligned}\)
        &
        $\mathcal Q^i$
        \\
        \bottomrule
    \end{tabular}}
\end{table*}

\subsection{Compressed \algname{PEARL} Variants}
Both \algname{PEARL-SGD} and \algname{PEARL-SVRG} communicate the full strategy increment at synchronization. 
We introduce \algname{C-PEARL-SGD} and \algname{C-PEARL-SVRG} as their compressed variants, where players communicate compressed local strategy increments rather than the full increments.

\subsubsection{\textbf{\algname{C-PEARL-SGD} algorithm}}  
\algname{C-PEARL-SGD} uses the same gradient estimator \eqref{eq:mini_batch_gradient} and local updates as \algname{PEARL-SGD}, but applies the compression maps \(\mathcal Q^i\) at synchronization. After the local updates, player \(i\) communicates the compressed strategy increment \(q_p^i=\mathcal Q^i(u_p^i)\). This can reduce communication overhead at the cost of additional error due to compression, which can enlarge the convergence neighborhood in the stochastic setting (Corollary~\ref{corollary:c-pearl-sgd-rate}).

\subsubsection{\textbf{\algname{C-PEARL-SVRG} algorithm}} 
\algname{C-PEARL-SVRG} is the compressed variant of \algname{PEARL-SVRG}, or equivalently, the variance-reduced variant of \algname{C-PEARL-SGD}. It is obtained from \cref{alg:unified-pearl} by choosing the variance-reduced gradient estimator \eqref{eq:pearl-svrg-gradient-estimator} and applying compression to the strategy increments. This algorithm combines the benefits of both variants by allowing compression to reduce communication overhead while retaining exact convergence to the equilibrium in the stochastic finite-sum setting, as shown in Corollary~\ref{corollary:c-pearl-svrg-rate}.

\section{Convergence Guarantees}
\label{section:convergence-guarantees}
In this section, we establish convergence guarantees for \algname{PEARL} variants. Algorithm~\ref{alg:unified-pearl} provides a generalized algorithmic skeleton that recovers all of them as special cases, and we follow the same viewpoint in our analysis. We identify a common bound~\eqref{eq:unified_round_descent} tracking the progress of the algorithms over communication rounds, and formalize it in Assumption~\ref{assm:key_assumption}. We then show in Lemma~\ref{lemma:methods-satisfy-key-assumption} that all four algorithms satisfy this bound with algorithm-specific constants. Our ``Master Theorem''~\ref{thm:main_theorem} uses this bound to establish the convergence guarantees for all variants.

\subsection{Unified Parametric Assumption}
We introduce a unified parametric assumption used in all our convergence results, similar to the unified approaches of \cite{gorbunov2020unified,gorbunov2021local,beznosikov2023stochastic, gower2019sgd}. Unlike these works, which bound the second moment of the stochastic gradient or operator at every iteration, we bound the distance to the equilibrium after each communication round. We do not require an auxiliary sequence of non-negative variables since we do not maintain separate memory vectors as in certain compression-based methods \cite{mishchenko2025distributed}.

\begin{assumption}
\label{assm:key_assumption}
Let \(\mathcal F_p\) be the $\sigma$-algebra generated by all randomness before communication round \(p\), so that \(\mathbf x_{\tau p}\) is \(\mathcal F_p\)-measurable, and let $\mathbb E_p[\cdot] :=\mathbb E[\cdot\mid \mathcal F_p].$ 
We assume that for a constant step-size $\gamma_k\equiv\gamma>0$, there exist $A_\tau,B_\tau\geq0$, possibly depending on the synchronization interval $\tau$ and $\gamma$, such that, for every communication round $p$, the given algorithm satisfies
\begin{align}
\notag
    \mathbb{E}_p\left[\|\mathbf{x}_{\tau(p+1)} - \mathbf{x}_\star\|^2\right] &\le \|\mathbf{x}_{\tau p} - \mathbf{x}_\star\|^2
    - 2\gamma \tau \left\langle F(\mathbf{x}_{\tau p}), \mathbf{x}_{\tau p} - \mathbf{x}_{\star}\right\rangle\\
    &\quad+ \gamma^2 \tau^2 \Bigl(A_\tau\left\langle F(\mathbf{x}_{\tau p}), \mathbf{x}_{\tau p} - \mathbf{x}_{\star}\right\rangle+B_\tau\Bigr).
    \label{eq:unified_round_descent}
\end{align}
\end{assumption}
Here, \(\mathbb E_p\) denotes expectation over the randomness generated during communication round \(p\), including the mini-batch sampling used to construct \(g_k^i\) and, for the compressed variants, the randomness of the compression maps \(\mathcal Q^i\). In \eqref{eq:unified_round_descent}, \(A_\tau\) affects the contraction factor and step-size condition, while \(B_\tau\) is the stochastic error term and satisfies \(B_\tau=0\) in the deterministic and variance-reduced settings. We now present how \eqref{eq:unified_round_descent} is satisfied for each of the proposed \algname{PEARL} methods.

\begin{lemma}
\label{lemma:methods-satisfy-key-assumption}
Assume \hyperref[assumption:convexity]{\textbf{\emph{(CVX)}}},
\hyperref[assumption:smoothness]{\textbf{\emph{(SM)}}},
\hyperref[assumption:quasi-strong-monotonicity]{\textbf{\emph{(QSM)}}},
and \hyperref[assumption:star-cocoercivity]{\textbf{\emph{(SCO)}}}
hold. Let \(\kappa=\ell/\mu\), \(\sigma^2=\sum_{i=1}^n\sigma_i^2\),
\begin{center}
\footnotesize
\(\displaystyle
\begin{aligned}
A_\tau^{\mathrm{GD}}
&:=\ell+\frac{2(\tau-1)}{\tau}L_{\max}\sqrt{\frac{\kappa}{3}},\\
B_\tau^{\mathrm{GD}} &:=\frac{\sigma^2}{\tau} \left[1+(\tau-1)L_{\max}\left(4\gamma\tau+\frac{1+4\gamma\tau^2L_{\max}}
{2\tau\sqrt{\ell\mu/3}}\right)\right],\\
U_\tau &:=\frac{\sigma^2}{\tau} \left(1+2(\tau-1)(\tau+1)\gamma L_{\max}\right),\\
A_\tau^{\mathrm{VR}} &:=4\ell+\frac{4(\tau-1)}{\tau} L_{\max}\sqrt{\frac{\kappa}{3}}.
\end{aligned}
\)
\end{center}
and $L_{\mathrm{max}} = \max\{L_1,\dots,L_n\}$. Then the following instances of Algorithm~\ref{alg:unified-pearl} satisfy
Assumption~\ref{assm:key_assumption} with:

\begin{itemize}[leftmargin=0.8em]

\item \textbf{\algname{PEARL-SGD}.}
\(A_\tau=A_\tau^{\mathrm{GD}}\) and
\(B_\tau=B_\tau^{\mathrm{GD}}\), provided that \hyperref[assumption:noise]{\textbf{\emph{(BV)}}} holds and
\(0<\gamma_k\equiv\gamma\leq
\frac{1}{\ell\tau+2(\tau-1)L_{\max}\sqrt{\kappa}}\).

\item \textbf{\algname{PEARL-SVRG}.}
\(A_\tau=A_\tau^{\mathrm{VR}}\) and \(B_\tau=0\), provided that Assumption~\ref{assumption:finite-sum-smoothness} holds and \(0<\gamma_k\equiv\gamma\leq\frac{1}{2L_{\max}\tau}\).

\item \textbf{\algname{C-PEARL-SGD}.}
\(A_\tau=A_\tau^{\mathrm{GD}}+\omega\ell\) and
\(B_\tau=B_\tau^{\mathrm{GD}}+\omega U_\tau\), provided that \hyperref[assumption:noise]{\textbf{\emph{(BV)}}} and
Assumption~\ref{assumption:unbiased-compression} hold, and
\(0<\gamma_k\equiv\gamma\leq
\frac{1}{\ell\tau+2(\tau-1)L_{\max}\sqrt{\kappa}}\).

\item \textbf{\algname{C-PEARL-SVRG}.}
\(A_\tau=A_\tau^{\mathrm{VR}}+4\omega\ell\) and \(B_\tau=0\), provided that Assumptions~\ref{assumption:finite-sum-smoothness}
and~\ref{assumption:unbiased-compression} hold, and
\(0<\gamma_k\equiv\gamma\leq\frac{1}{2L_{\max}\tau}\).
\end{itemize}
\end{lemma}

The differences among the algorithms are encoded in the constants \(A_\tau\) and \(B_\tau\). For the variance-reduced variants, \(B_\tau=0\), while compression introduces additional terms proportional to \(\omega\) in \(A_\tau\) for both compressed methods and in \(B_\tau\) for \algname{C-PEARL-SGD}.

\begin{table*}[t]
    \centering
    \caption{Overview of the \algname{PEARL} variants and their convergence guarantees with constant $\gamma$. ``Linear $+$ neighborhood'' denotes linear convergence to a neighborhood of the equilibrium, while ``Linear $+$ exact'' denotes linear convergence to the exact equilibrium. }
    \label{table:pearl-unified-results}
    \small
    \renewcommand{\arraystretch}{1.20}
    \setlength{\tabcolsep}{4.0pt}
    \begin{tabular}{
        >{\centering\arraybackslash}m{0.17\textwidth}
        >{\centering\arraybackslash}m{0.17\textwidth}
        >{\centering\arraybackslash}m{0.18\textwidth}
        >{\centering\arraybackslash}m{0.13\textwidth}
        >{\centering\arraybackslash}m{0.25\textwidth}
    }
        \toprule
        \textbf{Method}
        & \textbf{Result}
        & \textbf{Variance reduction}
        & \textbf{Compression}
        & \textbf{Convergence guarantee} \\
        \midrule

        \algname{PEARL-SGD}
        & \Cref{corollary:unified-pearl-sgd}
        & \(\times\)
        & \(\times\)
        & Linear \(+\) neighborhood \\

        \algname{PEARL-SVRG}
        & \Cref{corollary:pearl-svrg-rate}
        & \(\checkmark\)
        & \(\times\)
        & Linear \(+\) exact \\

        \algname{C-PEARL-SGD}
        & \Cref{corollary:c-pearl-sgd-rate}
        & \(\times\)
        & \(\checkmark\)
        & Linear \(+\) neighborhood \\

        \algname{C-PEARL-SVRG}
        & \Cref{corollary:c-pearl-svrg-rate}
        & \(\checkmark\)
        & \(\checkmark\)
        & Linear \(+\) exact \\

        \bottomrule
    \end{tabular}
\end{table*}

\vspace{-3mm}
\subsection{Convergence of \algname{PEARL} Variants}
We now present our general convergence result under Assumption~\ref{assm:key_assumption} and \hyperref[assumption:quasi-strong-monotonicity]{\textbf{\emph{(QSM)}}}. 

\begin{theorem}
\label{thm:main_theorem}
Assume \hyperref[assumption:quasi-strong-monotonicity]{\textbf{\emph{(QSM)}}} and Assumption~\ref{assm:key_assumption} hold. Suppose  $A_\tau\geq\mu$, $0<\gamma<\frac{2}{\tau A_\tau}$, and $\zeta:=2-\gamma\tau A_\tau >0$. Then, for every \(R\geq0\), the iterates generated by any specialization of Algorithm~\ref{alg:unified-pearl} satisfy
\begin{align*}
    \mathbb{E}\left[\|\mathbf{x}_{\tau R}-\mathbf{x}_\star\|^2\right]&\leq(1-\gamma\tau\mu\zeta)^R\|\mathbf{x}_{0}-\mathbf{x}_\star\|^2  +\frac{\gamma\tau B_{\tau}}{\mu \zeta}.
\end{align*}
\end{theorem}
This result gives linear convergence to a neighborhood of $\vx_\star$, and exact convergence when $B_\tau=0$. We now apply it with constants from  Lemma~\ref{lemma:methods-satisfy-key-assumption} to obtain convergence guarantees for each \algname{PEARL} variant.
\begin{corollary}
\label{corollary:unified-pearl-sgd}
Under the assumptions of Lemma~\ref{lemma:methods-satisfy-key-assumption}, if 
$0<\gamma<\min\left\{
\frac{1}{\ell\tau+2(\tau-1)L_{\max}\sqrt{\kappa}},
\frac{2}{\tau A_\tau^{\mathrm{GD}}}
\right\},$
and $\zeta^{\mathrm{GD}}:=2-\gamma\tau A_\tau^{\mathrm{GD}}$, then the iterates of \algname{PEARL-SGD} satisfy
\begin{align*}
    \mathbb E\left[
    \|\mathbf{x}_{\tau R}-\mathbf{x}_\star\|^2
    \right] &\le\left(1-\gamma\tau\mu\zeta^{\mathrm{GD}}\right)^R
    \|\mathbf{x}_0-\mathbf{x}_\star\|^2
    +\frac{\gamma\tau B_\tau^{\mathrm{GD}}}{\mu\zeta^{\mathrm{GD}}}.
\end{align*}
\end{corollary}
Corollary \ref{corollary:unified-pearl-sgd} gives linear convergence to a stochastic neighborhood and recovers the convergence guarantee of the original \algname{PEARL-SGD}~\cite[Theorem~3.4]{yoon2025multiplayer} (up to some typo corrections in constants). When $\tau=1$, the rate is consistent with the analysis of SGDA~\cite{beznosikov2023stochastic}, while increasing $\tau$ reduces the communication frequency, yielding the communication gain established in~\cite[Corollary~3.5]{yoon2025multiplayer}. 

\begin{corollary}
\label{corollary:pearl-svrg-rate}
Under the assumptions of Lemma~\ref{lemma:methods-satisfy-key-assumption}, if 
$0<\gamma<\min\left\{
\frac{1}{2L_{\max}\tau},
\frac{2}{\tau A_\tau^{\mathrm{VR}}}
\right\},$
and $\zeta^{\mathrm{VR}}:=2-\gamma\tau A_\tau^{\mathrm{VR}}$, then the iterates of \algname{PEARL-SVRG} satisfy
\begin{align*}
    \mathbb E\left[\left\|\mathbf{x}_{\tau R}-\mathbf{x}_\star\right\|^2
    \right]\leq\left(1-\gamma\tau\mu\zeta^{\mathrm{VR}}\right)^R
    \left\|\mathbf{x}_0-\mathbf{x}_\star\right\|^2.
\end{align*}
\end{corollary}
Unlike \algname{PEARL-SGD} with a constant step-size and nonzero stochastic gradient variance, the bound contains no neighborhood term. Hence, \algname{PEARL-SVRG} converges linearly to the exact equilibrium, consistent with standard variance-reduced methods like SVRG or L-SVRG~\cite{gorbunov2020unified, palaniappan2016stochastic, beznosikov2023stochastic}.

\begin{corollary}
\label{corollary:c-pearl-sgd-rate}
Under the assumptions of Lemma~\ref{lemma:methods-satisfy-key-assumption}, if 
$0<\gamma<\min\bigl\{\frac{1}{\ell\tau+2(\tau-1)L_{\max}\sqrt{\kappa}},
\allowbreak\frac{2}{\tau(A_\tau^{\mathrm{GD}}+\omega\ell)}\bigr\},$
and $\zeta^{{\mathrm{CGD}}}:= 2-\gamma\tau(A_\tau^{\mathrm{GD}}+\omega\ell),$ then the iterates of \algname{C-PEARL-SGD} satisfy
\begin{align*}
   \mathbb E\left[\left\|\mathbf{x}_{\tau R}-\mathbf{x}_\star\right\|^2\right]
    &\leq\left(1-\gamma\tau\mu\zeta^{\mathrm{CGD}}\right)^R\left\|\mathbf{x}_0-\mathbf{x}_\star\right\|^2+\frac{\gamma\tau  (B_\tau^{\mathrm{GD}}+\omega U_\tau)}{\mu\zeta^{\mathrm{CGD}}}.
\end{align*}
\end{corollary}
\algname{C-PEARL-SGD} converges linearly to a neighborhood of the equilibrium, with a potentially larger neighborhood bound than \algname{PEARL-SGD} due to the additional error from compression. When $\omega=0$, we recover the convergence guarantee of \algname{PEARL-SGD}. When deterministic gradients are used, both $B_\tau^{\mathrm{GD}}$ and $U_\tau$ vanish, and \algname{C-PEARL-SGD} converges linearly to the exact equilibrium despite compression. 

\begin{corollary}
\label{corollary:c-pearl-svrg-rate}
Under the assumptions of Lemma~\ref{lemma:methods-satisfy-key-assumption}, suppose
$0<\gamma<\min\bigl\{\frac{1}{2L_{\max}\tau}, \allowbreak\frac{2}{\tau(A_\tau^{\mathrm{VR}}+4\omega\ell)} \bigr\},$
and $\zeta^{\mathrm{CVR}}:=2-\gamma\tau(A_\tau^{\mathrm{VR}}+4\omega\ell).$ Then the iterates of \algname{C-PEARL-SVRG} satisfy
\begin{align*}
    \mathbb E\left[\left\|\mathbf{x}_{\tau R}-\mathbf{x}_\star\right\|^2\right]
    \leq \left(1-\gamma\tau\mu\zeta^{\mathrm{CVR}}\right)^R
    \left\|\mathbf{x}_0-\mathbf{x}_\star\right\|^2.
\end{align*}
\end{corollary}
\algname{C-PEARL-SVRG} converges linearly to the equilibrium despite using stochastic gradients and communicating compressed strategy increments, and when \(\omega=0\), we recover the convergence guarantee of \algname{PEARL-SVRG}.

\section{Numerical Experiments}
\label{section:experiments}
In this section, we present experiments to illustrate our theoretical results and demonstrate the practical benefits of our proposed algorithms. We consider two settings: a quadratic $n$-player game and a distributed spectrum sensing problem in cognitive radio networks.
% \vspace{-2mm}
\subsection{Quadratic $n$-player game} 
We consider a quadratic $n$-player game, which models a broad range of problems considered in game theory literature  \cite{scutari2010convex, chen2014multitask, nassif2016multitask, YeHu2017_distributed, SalehisadaghianiPavel2016_distributed}. Player $i \in [n]$ has objective 
\begin{equation}
\label{exp:quadnPlayerGame}
    f_i(x^i; x^{-i}):=\frac{1}{M}\sum_{m=1}^{M} f_{i,m}(x^i;x^{-i}),
\end{equation}
where $f_{i,m}(x^i;x^{-i}) =\frac12\langle x^i,\mathbf A_{i,m}x^i\rangle+\sum_{\substack{1\le j\le n\\ j\neq i}}\langle x^i,\mathbf B_{i,j,m}x^j\rangle+\langle a_{i,m},x^i\rangle$ with $\vA_{i,m} \in \reals^{d_i \times d_i}, \mathbf B_{i,j,m}\in\mathbb R^{d_i\times d_j}$ and $a_{i,m}\in\mathbb R^{d_i}$. We take $n=5$, $d_i = d = 10$, and $M=100$. The matrices $\vA_{i,m}$ and $\vB_{i,j,m}$ (for $i<j$)  are random symmetric $d\times d$ matrices; we choose their eigenvalues uniformly randomly in $[\mu_A, L_A]$ and $[\mu_B, L_B]$, respectively, and their orthogonal eigenvectors by taking the QR-decomposition of a $d\times d$ random matrix with i.i.d.\ standard normal entries.
We set $\vB_{j,i,m} = -\vB_{i,j,m}^{\mathsf{T}}$ for $j>i$ to enforce strong monotonicity of $F$, and take $\mu_A = 0.01$, $L_A = 0.1$, $\mu_B = 0$, and $L_B = 1.0$.
We use \texttt{Rand-$k$} compression \cite{beznosikov2023stochastic} with $k=5$ for compressed algorithms. \texttt{Rand-$k$} selects \(k\) coordinates uniformly without replacement, scales the selected values by \(d/k\), and sets the remaining coordinates to zero.

\begin{figure*}[h]
\centering
\subfloat[Error vs.\ Communication]{%
\includegraphics[width=0.31\textwidth]{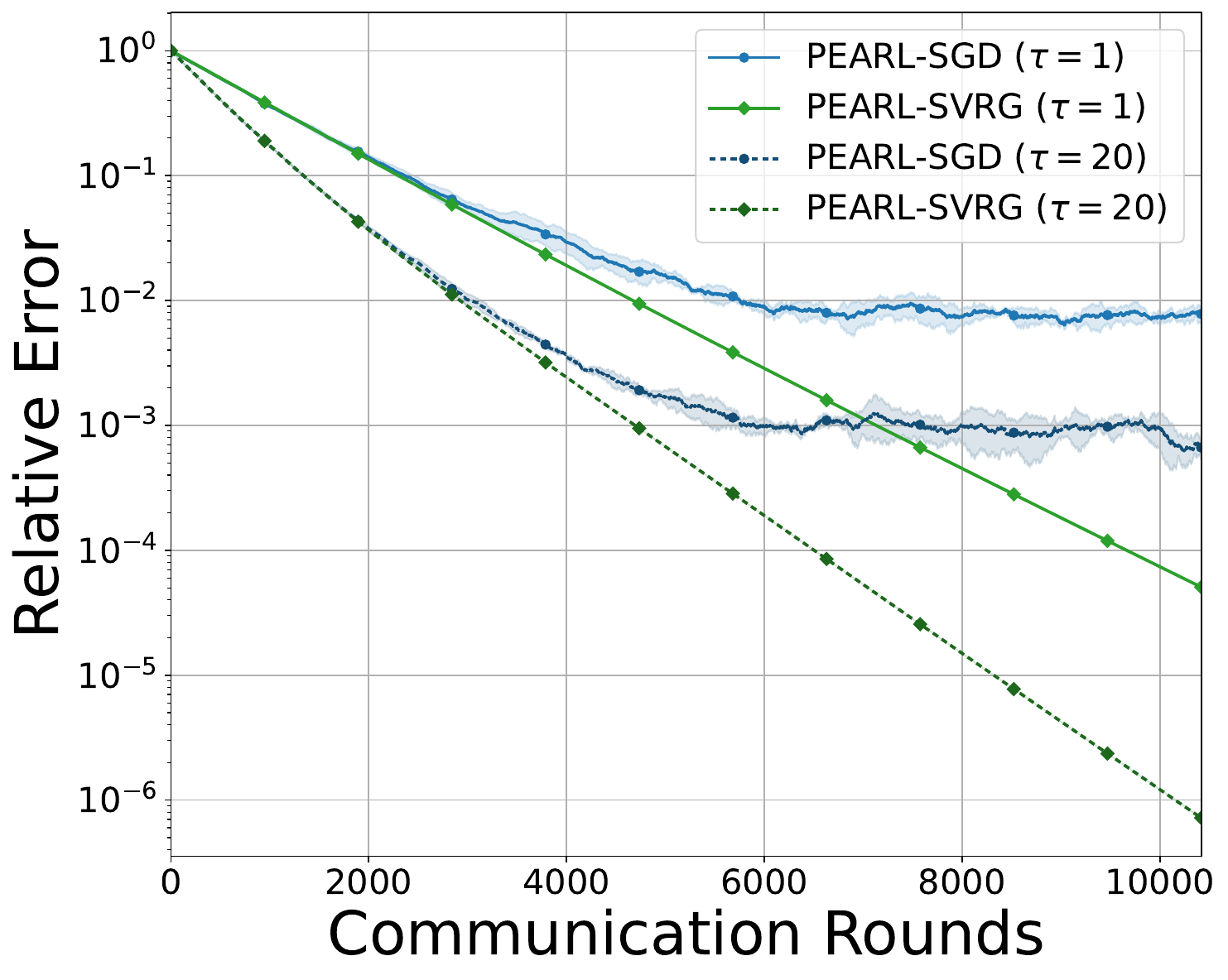}%
\label{fig:np-variance-rounds-comparison}}
\hfil
\subfloat[Error vs.\ Oracle calls]{%
\includegraphics[width=0.31\textwidth]{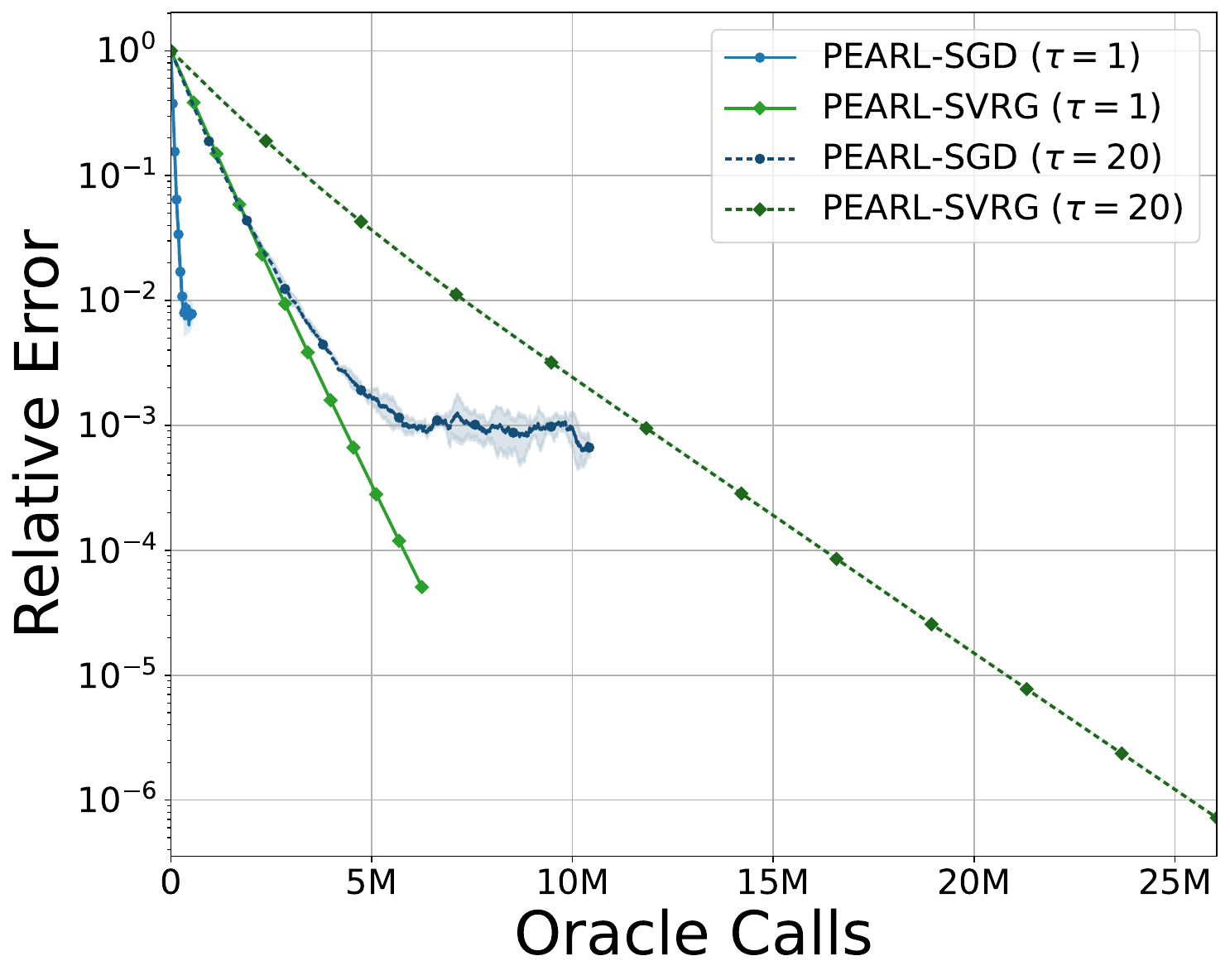}%
\label{fig:np-variance-oracle-comparison}}
\hfil
\subfloat[Error vs.\ Communicated bits]{%
\includegraphics[width=0.31\textwidth]{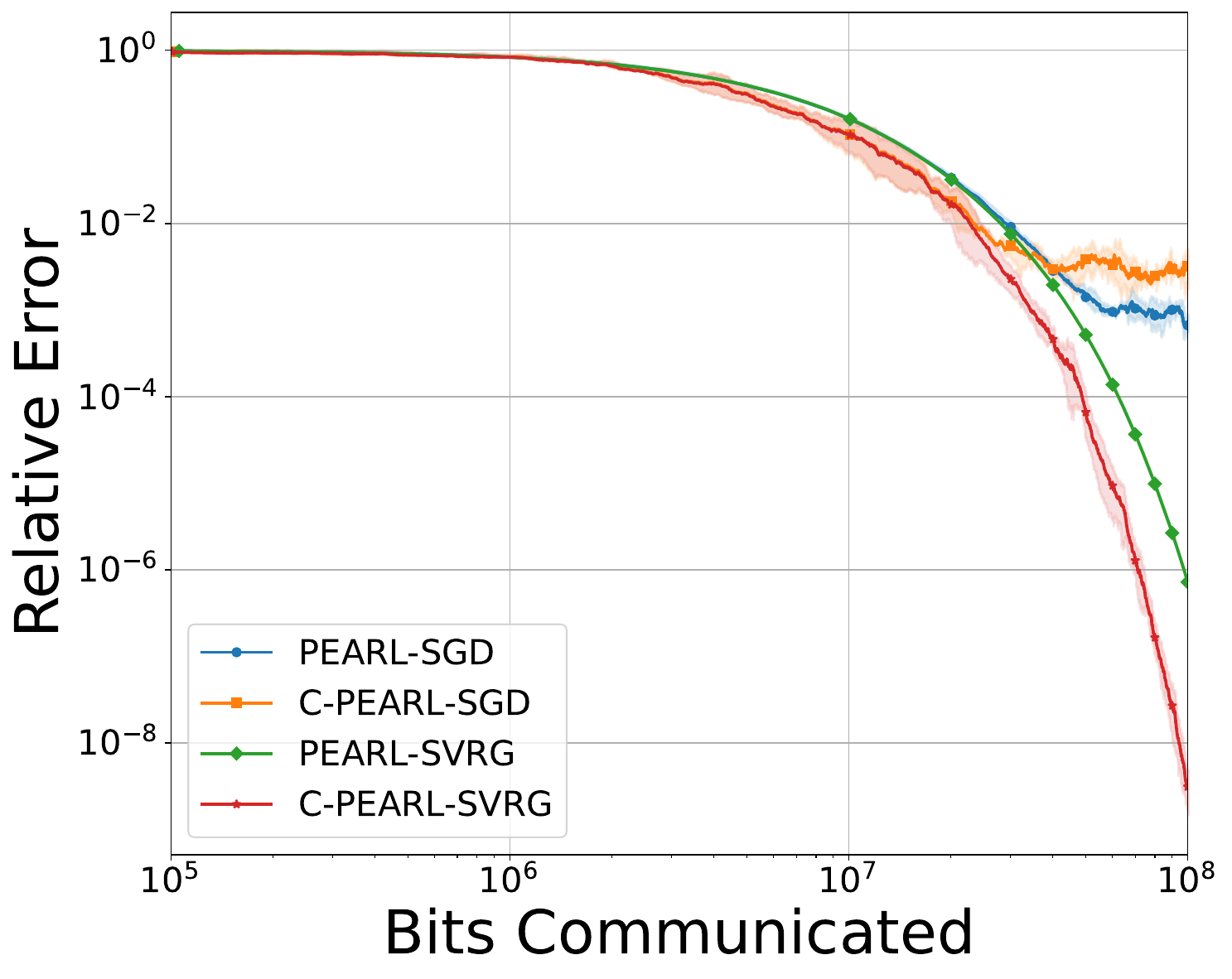}%
\label{fig:np-bit-comparison}}
\caption{Performance plot of \algname{PEARL} algorithm variants on the quadratic \(n\)-player game~\eqref{exp:quadnPlayerGame} using tuned step-sizes. The first two plots illustrate the convergence of \algname{PEARL-SVRG} to the exact solution, traded off with increased oracle complexity. The rightmost plot compares all variants with $\tau=20$, illustrating that algorithms using compression can converge with less bit communication budget compared to their non-compressed counterparts.
}
\label{fig:np-game-results}
\end{figure*}

We measure the performance of the \algname{PEARL} variants using the relative error $\frac{\sqnorm{\vx_{\tau p} - \vx_\star}}{\sqnorm{\vx_0 - \vx_\star}}$ and run all algorithms with a fixed communication budget of $10^8$ bits. For each algorithm with $\tau$ fixed, we select the best step-size from $\{10^{-2}, 10^{-3}, 10^{-4}, 10^{-5}\}$ based on the mean final error over three random trials, and plot the mean relative error with standard deviation (shaded regions). \Cref{fig:np-variance-rounds-comparison} demonstrates that \algname{PEARL-SGD} converges to a neighborhood, while \algname{PEARL-SVRG} converges linearly to the exact solution. For both algorithms, larger $\tau$ enables faster convergence for the same number of communication rounds. However, \algname{PEARL-SVRG} uses more oracle calls per round since it computes a full-batch gradient at the beginning of each round, as shown in \cref{fig:np-variance-oracle-comparison}.

When we additionally compare the compressed variants, \cref{fig:np-bit-comparison} shows that they require fewer communicated bits than their non-compressed counterparts during the initial phase where the error decreases nearly linearly. \algname{C-PEARL-SGD} converges to a larger neighborhood than \algname{PEARL-SGD} and has the worst final performance, whereas \algname{C-PEARL-SVRG} converges to the solution and performs best at high accuracy in terms of communicated bits.

\paragraph{Remark.} \cref{fig:compression-scheme-comparison} uses the same setup as in this experiment with $\tau=5$, $\gamma=10^{-3}$, $L_B=2$ and $5\times 10^4$ communication rounds, and uses full-batch gradients without stochasticity.

\begin{figure*}
\centering
\subfloat[Error vs.\ Communication]{%
\includegraphics[width=0.31\textwidth]{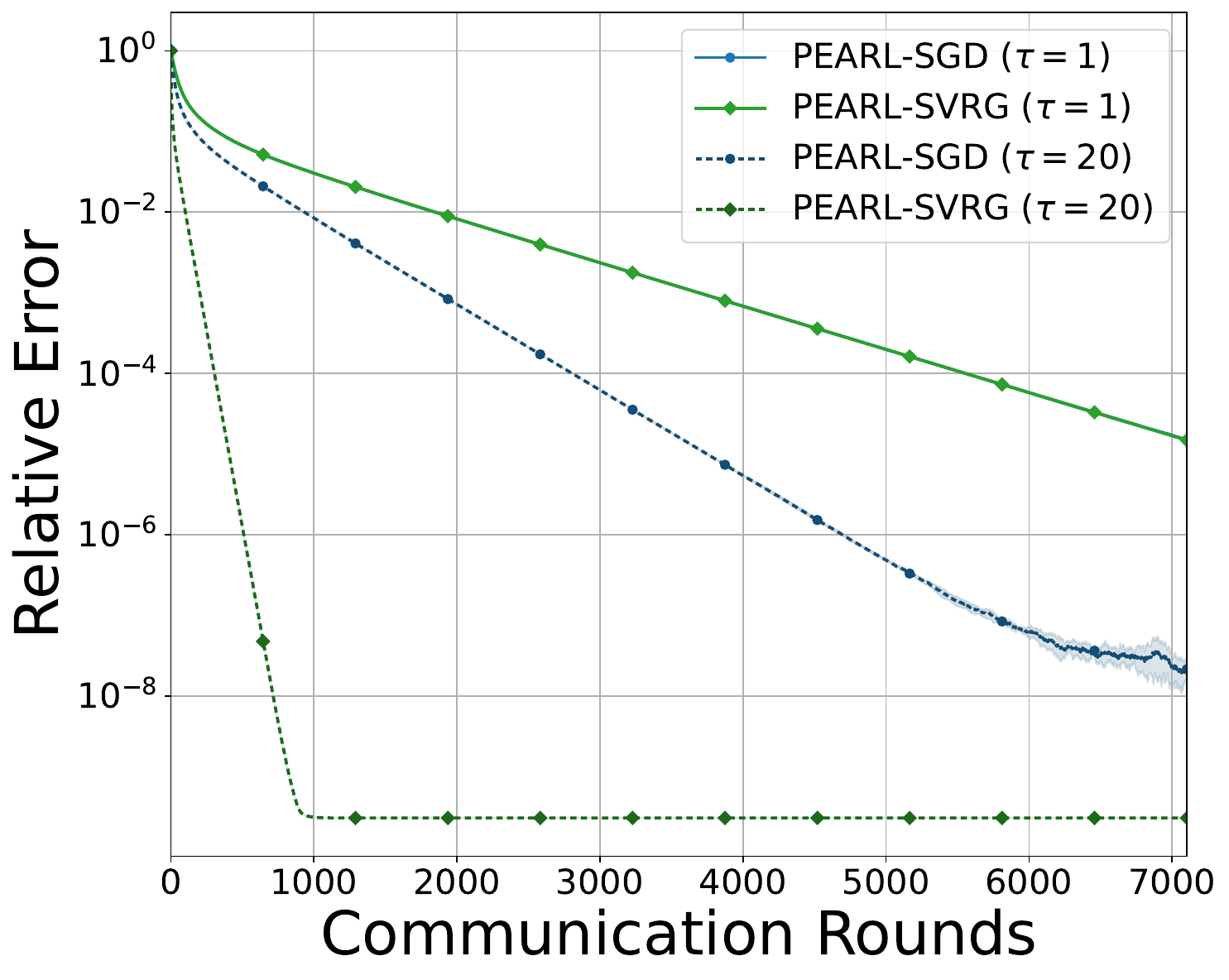}%
\label{fig:ss-variance-rounds-comparison}}
\hfil
\subfloat[Error vs.\ Oracle calls]{%
\includegraphics[width=0.31\textwidth]{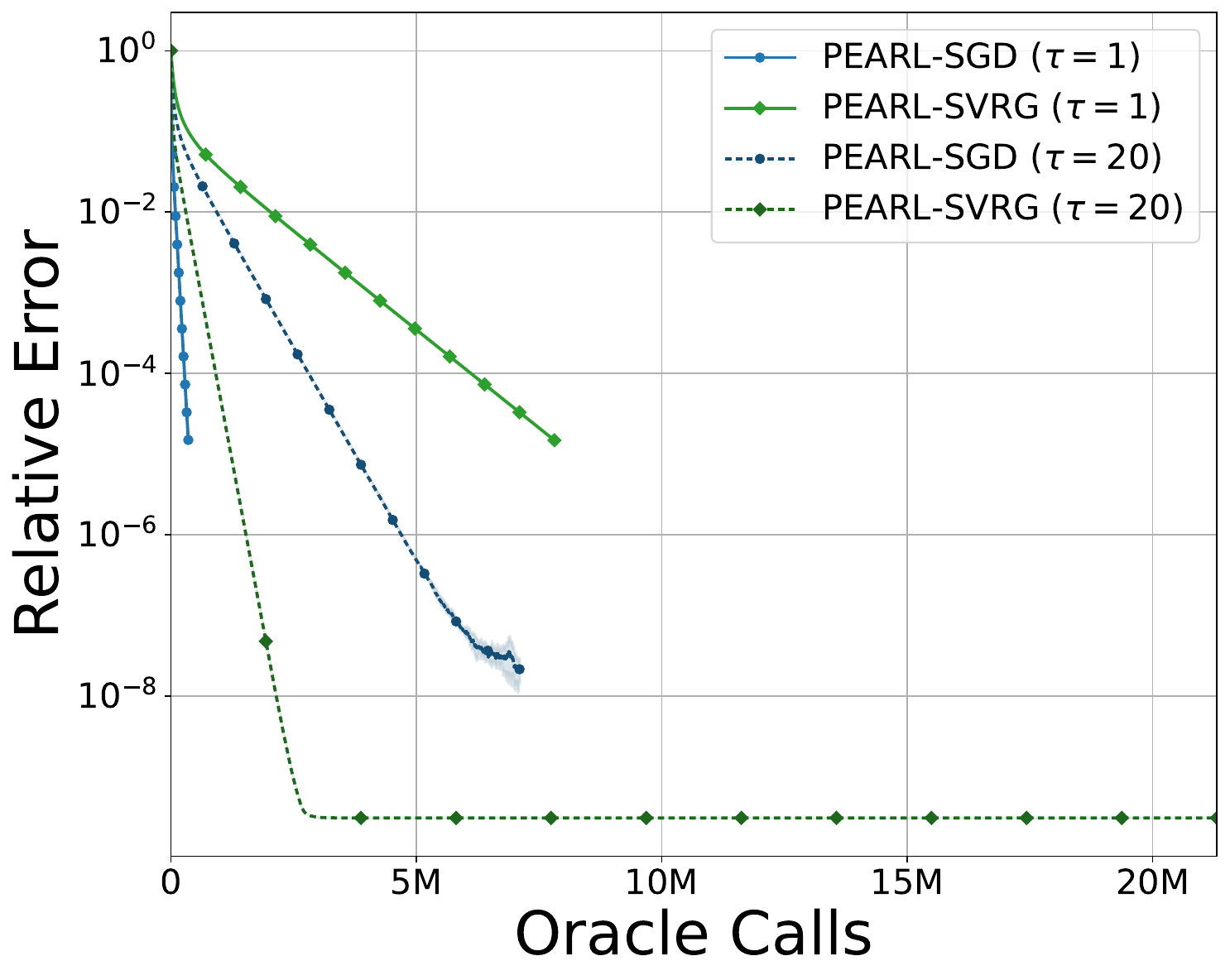}%
\label{fig:ss-variance-oracle-comparison}}
\hfil
\subfloat[Error vs.\ Communicated bits]{%
\includegraphics[width=0.31\textwidth]{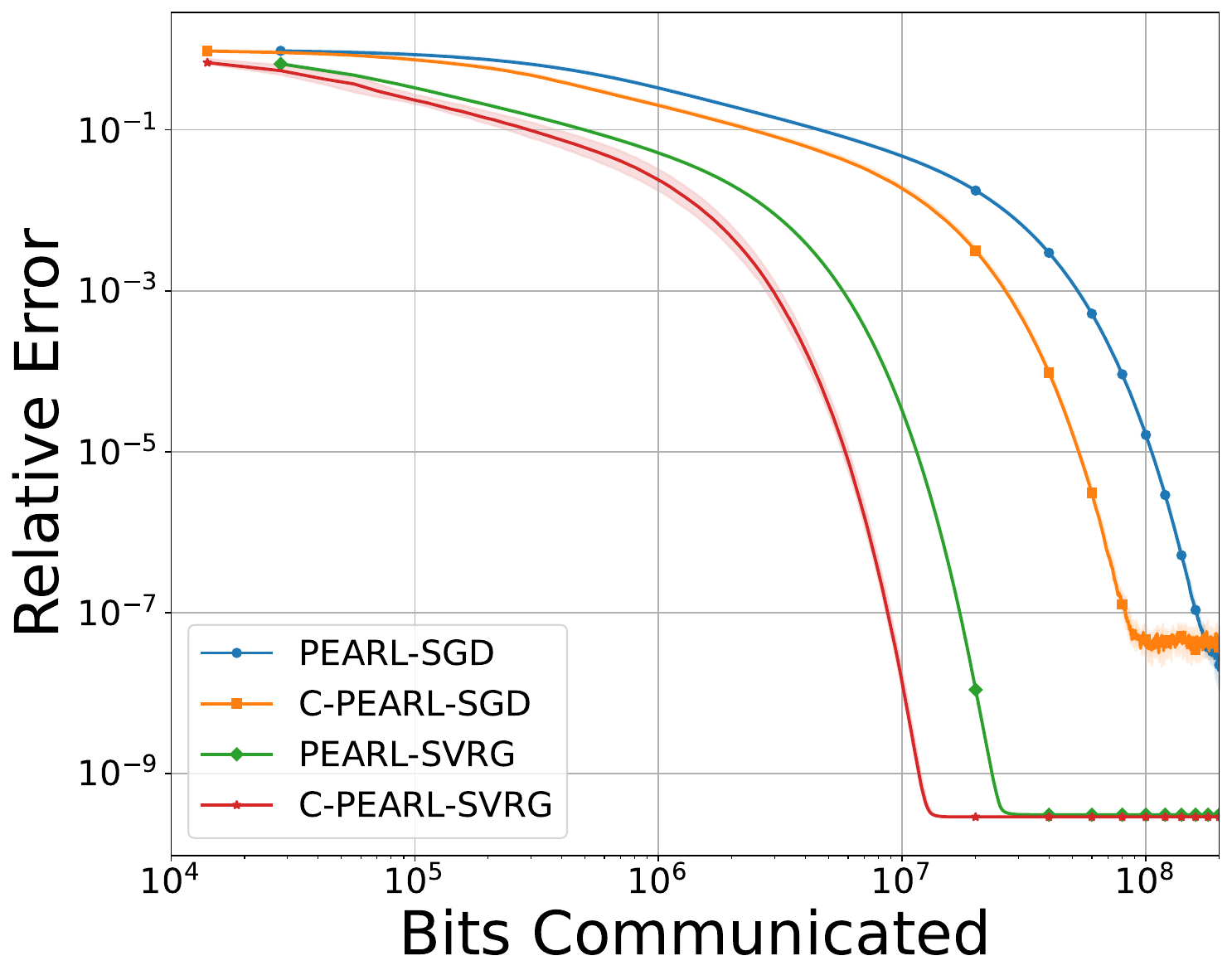}%
\label{fig:ss-bit-comparison}}
\caption{Performance plot of \algname{PEARL} algorithm variants on the distributed spectrum sensing problem~\eqref{eq:spectrum-sensing-game-objective} using tuned step-sizes. The first two plots illustrate the convergence of \algname{PEARL-SVRG} to the exact solution, traded off with increased oracle complexity. The rightmost plot compares all variants with $\tau=20$, illustrating that algorithms using compression can converge with less bit communication budget compared to their non-compressed counterparts.
}
\label{fig:spectrum-sensing-results}
\end{figure*}
\vspace{-2mm}
\subsection{Distributed spectrum sensing}
We consider the distributed spectrum sensing problem of~\cite{chen2014multitask}. We are given primary users \(q = 1,\ldots,N_P\) that are licensed to transmit signals over certain frequency bands. They are fixed environmental factors, not the players of the game. The players are the secondary users \(i = 1,\ldots,N_S\), who try to estimate primary users' spectra to identify idle bands. The power spectrum of primary user \(q\) is represented as $S_q(f)=\sum_{m=1}^{N_B}\alpha_{q,m}\phi_m(f)$, where each \(\phi_m\) is a Gaussian basis function, $f$ is the frequency variable and \(\alpha_{q,m}\) is the corresponding coefficient. The secondary users aim to accurately estimate the ground-truth coefficient vector of all primary users \(\alpha_\star\in\mathbb{R}^{N_PN_B}\), with each secondary user's estimate denoted by $x^i\in\mathbb{R}^{N_PN_B}$. Each secondary user has \(N_R\) antennas. Each antenna $\ell$ measures $r_{i\ell,s} =\Phi_{i\ell,s}\alpha_\star+z_{i\ell,s} \in \reals^{N_F}$ for sample $s\in[M]$, where $z_{i\ell,s}$ is the measurement noise, $\Phi_{i\ell,s}\in\mathbb{R}^{N_F\times N_PN_B}$ is the sensing matrix modeling the attenuation of primary-user signals as they propagate to the secondary user, and $\widehat{\Phi}_{i\ell,s}\in\mathbb{R}^{N_F\times N_PN_B}$ is its estimate.
 The objective function of the $i$-th player (secondary user) is given by 
 \begin{align}
\begin{aligned}
    f_i(x^i;x^{-i})&=\underbrace{\frac{1}{M}\sum_{s=1}^{M}\sum_{\ell=1}^{N_R} \frac12\left\|r_{i\ell,s}-\widehat{\Phi}_{i\ell,s}x^i\right\|^2}_{=:J_{i,1}(x^i)} + \underbrace{\frac{\eta N_R}{2}\sum_{j\in\mathcal N_i}\rho_{ij}\left\|x^i-x^j\right\|^2}_{=:J_{i,2}(x^i;x^{-i})}
\end{aligned}
\label{eq:spectrum-sensing-game-objective}
\end{align}
where \(J_{i,1}(x^i)\) is the local spectrum estimation error and the \(J_{i,2}(x^i;x^{-i})\) term induces similarity between neighboring users' estimates. Here, \(\mathcal N_i\) is the set of neighbors of player \(i\), \(\rho_{ij}\geq0\) is the neighbor weight, and \(\eta>0\) controls cooperation.

We take \(N_S=10\), \(N_P=2\), \(N_R=10\), \(N_B=4\), \(N_F=20\), and \(M=100\), so each strategy has dimension \(d=N_PN_B=8\). We set
\(\alpha_\star^{(1)}=(0.3,0.3,0.3,0)\), \(\alpha_\star^{(2)}=(0,0.3,0.3,0.3)\), and \(\alpha_\star=(\alpha_\star^{(1)},\alpha_\star^{(2)})\). The Gaussian basis functions, frequency grid, and sensing matrices are generated from hypothetical 2D locations using a nominal path-loss model with Gaussian noise. We determine $\cN_i$ by the ring graph, where \(\rho_{i,i-1}=\rho_{i,i+1}=1/2\) and \(\rho_{ij}=0\) otherwise. We take \(\eta=1.0\), use \texttt{Rand-\(k\)} compression with \(k=4\) and mini-batch size \(5\), and plot the mean of the relative error $\frac{\sqnorm{\vx_{\tau p} - \vx_\star}}{\sqnorm{\vx_0 - \vx_\star}}$ with standard deviation (shaded region) over five trials for \(\tau\in\{1,20\}\). For each case, we choose the best step-size from \(\{10^{-3},10^{-4},10^{-5},10^{-6}\}\) under a fixed communication budget of \(2\times10^8\) bits.

We observe in \cref{fig:ss-variance-rounds-comparison} that \algname{PEARL-SVRG} with $\tau=20$ converges linearly to the solution in fewer than \(10^3\) communication rounds, while other algorithms converge more slowly, and \algname{PEARL-SGD} with $\tau=20$ eventually stagnates in a neighborhood. Although variance reduction requires more oracle calls per round, \cref{fig:ss-variance-oracle-comparison} shows that \algname{PEARL-SVRG} still reaches the exact solution with fewer oracle calls in this experiment. Finally, \cref{fig:ss-bit-comparison} shows that the compressed variants of \algname{PEARL-SGD} and \algname{PEARL-SVRG} reach the solution or its neighborhood using fewer communicated bits, as expected.

\section{Conclusion}
We develop generalized and advanced \algname{PEARL} methods for MpFL based on an equivalent increment-sharing reformulation of \algname{PEARL-SGD}, incorporating compressed communication and variance reduction. We provide a unified convergence analysis for all variants, including the newly developed \algname{C-PEARL-SGD}, \algname{PEARL-SVRG} and \algname{C-PEARL-SVRG}. Numerical experiments support the theory and illustrate the practical benefits of the proposed methods. Our results offer several natural extensions for future work, including different variance-reduced estimators or baseline algorithms using history-dependent update directions such as momentum or anchoring methods.

\subsection*{Acknowledgments}
The authors’ contribution to this work was supported by NSF CCF 2504626, NSF CAREER
2542902, and AFOSR YIP FA9550-26-1-B353.

\subsection*{AI use statement}
In this work, we used generative AI tools to assist with polishing the paper, including improving clarity and LaTeX formatting. We also used generative AI tools to assist with mathematical derivations and proof checking, as well as the development and debugging of experimental code. We did not use generative AI tools to generate the ideas or concepts of the project, or to generate citations. All AI-assisted mathematical arguments and proofs were checked, and all AI-assisted code was reviewed and tested by the authors. Generative AI was not used to generate experimental data or results. We take responsibility for the final content of this work, including all text, claims, and artifacts produced with the aid of generative AI.

\bibliographystyle{plain}
\bibliography{references}
\newpage
\appendix
\section*{Appendix}
We organize the appendix as follows: Section \ref{appendix:related_works} provides a discussion on related works, Section \ref{appendix:convergence-guarantees} provides the proofs of the main convergence results in Section \ref{section:convergence-guarantees}, and Section \ref{appendix:details_exact_convergence_increment_sharing} provides details on the increment-sharing mechanism as well as its comparison with algorithms in distributed minimization discussed in Section \ref{section:MpFL-concept}.

\tableofcontents

\section{Related Works}
\label{appendix:related_works}

\subsection{Multiplayer Games under Distributed Systems}

\paragraph{Multiplayer Federated Learning (MpFL).}
Multiplayer games have been studied in different distributed settings, including MpFL~\cite{yoon2025multiplayer}, which is the focus of this paper. A main objective of MpFL is to find an equilibrium in as few communication rounds as possible. To reduce communication, \algname{PEARL-SGD}~\cite{yoon2025multiplayer} allows players to perform several local updates between communication rounds. However, excessive local updates may result in \emph{player drift}, where players' actions approach minimizers of individual objectives while fixing the outdated strategies of the other players, which causes the convergence to equilibrium to fail. This issue is addressed in \algname{PEARL-Prox}~\cite{yoon2025pearl} by incorporating proximal updates into each player's local training loop. A closely related work, Decoupled SGD~\cite{zindari2025decoupled}, also uses the idea of performing several local updates while keeping the strategies of the other players fixed, and exhibits the same update rule as \algname{PEARL-SGD}. However, unlike \algname{PEARL-SGD}, its analysis focuses on communication efficiency gain for weakly coupled games in the deterministic setting. More recently, \cite{luo2026efficient} extended this line of work by developing decoupled methods for distributed saddle-point problems and monotone variational inequalities.

\paragraph{Distributed Nash Equilibrium Seeking.}
Another related setting is distributed Nash equilibrium seeking~\cite{kannan2012distributed,lei2020synchronous}, where each player performs computations using its own objective while exchanging partial information of action vectors with other players to find an equilibrium. While MpFL focuses on reducing communication between players under the orchestration of a central server, works within this literature, including~\cite{TatarenkoShiNedic2021_geometric,nguyen2025distributed}, generally study communication over decentralized networks, with greater emphasis on convergence under different network topologies rather than communication complexity.

\subsection{Overview of Compression Techniques in Federated Learning and Multiplayer Games}

\paragraph{Compression in Distributed and Federated Optimization.}
The technique of compressing the communicated information 
has been extensively studied in distributed and federated optimization. Common approaches include stochastic quantization~\cite{alistarh2017qsgd,seide20141} and sparsification~\cite{stich2018sparsified,lin2018deep}, along with error-feedback mechanisms~\cite{karimireddy2019error,condat2022ef} to address compression errors. Compression has also been incorporated into federated learning algorithms that perform multiple local updates between communication rounds~\cite{reisizadeh2020fedpaq,haddadpour2021federated}. These methods reduce communication overhead by transmitting compressed gradients/model updates instead of full vectors. Beyond the choice of compression schemes, another important consideration is the choice of quantities to compress. 
For example, in distributed minimization, directly compressing local gradients can prevent exact convergence with a constant step-size, even when deterministic gradients are used, since gradients of local objectives need not vanish at the global minimizer. DIANA~\cite{mishchenko2025distributed} addresses this issue by maintaining auxiliary reference sequences and compressing the differences between local gradients and these references, so that the compression error vanishes as the algorithm converges. We provide a more detailed explanation of this in Section~\ref{appendix:details_exact_convergence_increment_sharing}.

\paragraph{Compression in Nash Equilibrium Seeking.}
Compression has also been explored in distributed Nash equilibrium seeking. One such method is C-DNES~\cite{chen2024linear}, which, like DIANA, maintains an auxiliary reference sequence and compresses the difference between each player's estimate of the joint strategy vector and this reference. On the other hand, in our increment-sharing reformulation of \algname{PEARL-SGD}, an additional reference sequence is not needed; see Appendix~\ref{appendix:details_exact_convergence_increment_sharing} for explanation. 

\section{Missing proofs from Section \ref{section:convergence-guarantees}}
\label{appendix:convergence-guarantees}
\subsection{Auxiliary Bounds for \algname{PEARL-SGD}}

Here we present Lemmas \ref{lemma:proof-outline-basic-bound} and \ref{lemma:local-GDA-iterate-difference-bound-stochastic-case} from \cite{yoon2025multiplayer}, which were used in establishing convergence guarantees for \algname{PEARL-SGD}. Lemma \ref{lemma:proof-outline-basic-bound} shows that one communication round of \algname{PEARL-SGD} behaves like a large single descent step. However, it also introduces local error terms caused by running multiple SGD steps locally, and Lemma~\ref{lemma:local-GDA-iterate-difference-bound-stochastic-case} bounds these local terms. 
 
 \begin{lemma}
\label{lemma:proof-outline-basic-bound}
\cite[Lemma~C.1]{yoon2025multiplayer} Assume \hyperref[assumption:smoothness]{\textbf{\emph{(SM)}}}, and let $L_{\mathrm{max}} = \max\{L_1,\dots,L_n\}$.
Let $0 \le p \le R-1$ and suppose $\gamma_k \equiv \gamma > 0$ for $k=\tau p, \dots, \tau(p+1) - 1$.
Then for any $\alpha > 0$, we have
\begin{align*}
\textstyle
    \mathbb{E}_p \left[\sqnorm{\vx_{\tau (p+1)} - \vx_\star} \right] 
    &\le (1 + (\tau -1)\alpha\gamma) \sqnorm{\vx_{\tau p} - \vx_\star} - 2\gamma\tau \inprod{F(\vx_{\tau p})}{\vx_{\tau p} - \vx_\star} \\
    & + \frac{\gamma L_{\mathrm{max}}^2}{\alpha} \sum_{j=\tau p+1}^{\tau p+\tau-1} \mathbb{E}_p \left[ \sqnorm{\vx_{\tau p} - \vx_{j}} \right]+ \mathbb{E}_p \left[ \sqnorm{\vx_{\tau p} - \vx_{\tau (p+1)}} \right] .
\end{align*}
\end{lemma}

\begin{lemma}
\label{lemma:local-GDA-iterate-difference-bound-stochastic-case}
\cite[Lemma~C.2]{yoon2025multiplayer} Assume~\hyperref[assumption:convexity]{\textbf{\emph{(CVX)}}}, \hyperref[assumption:smoothness]{\textbf{\emph{(SM)}}} and  \hyperref[assumption:noise]{\textbf{\emph{(BV)}}}.
For $i\in [n]$ and $0\le p \le R-1$, suppose $\gamma_k \equiv \gamma$ for $k=\tau p, \dots, \tau(p+1) - 1$, where $0 < \gamma \le \frac{1}{L_i} \min\left\{ 1, \frac{1}{\tau-1} \right\}$.
Then for $t=0,\dots,\tau$,
\begin{align*}
    &\mathbb{E}_p \left[ \sqnorm{x_{\tau p}^i - x_{\tau p + t}^i} \right] \le \gamma^2 t^2 \sqnorm{\nabla f_i(x_{\tau p}^i; x_{\tau p}^{-i})}+ \gamma^2 t \left(1 + 2(t-1)(t+1) \gamma L_i \right) \sigma_i^2 .
\end{align*}
\end{lemma}
Note: We use the convention $1/(\tau-1)=+\infty$ when $\tau=1$, so that $\min\{1,1/(\tau-1)\}=1$ in this case.

Lemma \ref{lemma:pearl-sgd-one-round-recursion} follows from Lemma \ref{lemma:proof-outline-basic-bound} and Lemma \ref{lemma:local-GDA-iterate-difference-bound-stochastic-case}. We use Lemma~\ref{lemma:pearl-sgd-one-round-recursion} throughout our analysis as it helps us verify \eqref{eq:unified_round_descent} from Assumption~\ref{assm:key_assumption} in the analysis of our special cases. 

\begin{lemma}
\label{lemma:pearl-sgd-one-round-recursion}
Suppose the assumptions of Lemma~\ref{lemma:proof-outline-basic-bound} and Lemma~\ref{lemma:local-GDA-iterate-difference-bound-stochastic-case} hold. Assume \hyperref[assumption:quasi-strong-monotonicity]{\textbf{\emph{(QSM)}}} and \hyperref[assumption:star-cocoercivity]{\textbf{\emph{(SCO)}}}. Let \(0\le p\le R-1\) be a fixed round index in \algname{PEARL-SGD}, and suppose
\(\gamma_k\equiv\gamma
    \le
    \frac{1}{\ell\tau+2(\tau-1)L_{\max}\sqrt{\kappa}}\) for
\(k=\tau p,\ldots,\tau(p+1)-1\). Then
\begin{align*}
\mathbb E_p \left[\left\|\mathbf{x}_{\tau(p+1)}-\mathbf{x}_\star\right\|^2 \right]&\le\left\|\mathbf{x}_{\tau p}-\mathbf{x}_\star\right\|^2-
2\gamma \tau\left\langle F(\mathbf{x}_{\tau p}),\mathbf{x}_{\tau p}-\mathbf{x}_\star\right\rangle\\
&+\gamma^2 \tau^2\left(\ell+\frac{2(\tau-1)}{\tau}
L_{\max}\sqrt{\frac{\kappa}{3}}\right)\left\langle F(\mathbf{x}_{\tau p}),\mathbf{x}_{\tau p}-\mathbf{x}_\star\right\rangle
\notag\\
&+\gamma^2\sigma^2 \tau(\tau-1)\left(
4\gamma\tau L_{\max}+\frac{L_{\max}}{2\tau\sqrt{\ell\mu/3}}
+\frac{2\gamma\tau L_{\max}^2}{\sqrt{\ell\mu/3}}\right)+\gamma^2\sigma^2 \tau,
\end{align*}
where \(\kappa=\ell/\mu\) and \(\sigma^2=\sum_{i=1}^n\sigma_i^2\).
\end{lemma}

\begin{proof}
For $\tau \geq 2$, the step-size condition of Lemma~\ref{lemma:local-GDA-iterate-difference-bound-stochastic-case} is satisfied since
$\gamma\le\frac{1}{\ell\tau+2(\tau-1)L_{\max}\sqrt{\kappa}}\le\frac{1}{L_{\max}}\min\left\{1,\frac{1}{\tau-1}\right\}.$ For $\tau=1$, the bound in Lemma~\ref{lemma:local-GDA-iterate-difference-bound-stochastic-case} holds for every $\gamma>0$ by \hyperref[assumption:noise]{\textbf{\emph{(BV)}}}.
Combining Lemmas~\ref{lemma:proof-outline-basic-bound} and
\ref{lemma:local-GDA-iterate-difference-bound-stochastic-case}, we obtain
\begin{align*}
\mathbb E\left[\left\|\mathbf{x}_{\tau(p+1)}-\mathbf{x}_\star\right\|^2\,\middle|\,\cF_p\right]&\leq\left(1+\alpha\gamma(\tau-1)\right)\left\|\mathbf{x}_{\tau p}-\mathbf{x}_\star\right\|^2-2\gamma\tau\left\langle
F(\mathbf{x}_{\tau p}),\mathbf{x}_{\tau p}-\mathbf{x}_\star\right\rangle\\
&\quad+\left(\gamma^2\tau^2+\frac{\gamma^3L_{\max}^2\tau^2(\tau-1)}{3\alpha}
\right)\left\|F(\mathbf{x}_{\tau p})\right\|^2\\
&\quad+\gamma^2\sigma^2\tau (\tau-1)\gamma L_{\max} \left(2(\tau+1)+\frac{L_{\max}}{2\alpha}+\frac{\gamma L_{\max}^2}{2\alpha}(\tau+1)^2\right)\\
&\quad+\gamma^2\sigma^2\tau\\
&\hspace{-1.6em}\overset{\hyperref[assumption:quasi-strong-monotonicity]{\textbf{\emph{(QSM)}}}, \hyperref[assumption:star-cocoercivity]{\textbf{\emph{(SCO)}}}}{\leq}\left\|\mathbf{x}_{\tau p}-\mathbf{x}_\star
\right\|^2-2\gamma \tau\left\langle F(\mathbf{x}_{\tau p}),\mathbf{x}_{\tau p}-\mathbf{x}_\star\right\rangle\\
&\quad+\left(\frac{\alpha\gamma(\tau-1)}{\mu}+\ell\gamma^2\tau^2+
\frac{\ell\gamma^3L_{\max}^2\tau^2(\tau-1)}{3\alpha}\right)\left\langle F(\mathbf{x}_{\tau p}),\mathbf{x}_{\tau p}-\mathbf{x}_\star
\right\rangle\\
&\quad+\gamma^2\sigma^2\tau(\tau-1)\gamma L_{\max}\left(2(\tau+1)+\frac{L_{\max}}{2\alpha}+\frac{\gamma L_{\max}^2}{2\alpha}(\tau+1)^2\right)\\
&\quad+\gamma^2\sigma^2\tau.
\end{align*}
Now substituting \( \alpha =\gamma\tau L_{\max}\sqrt{\frac{\ell\mu}{3}}\) in the above inequality,
\begin{align*}
\mathbb E_p \left[\left\|\mathbf{x}_{\tau(p+1)}-\mathbf{x}_\star
\right\|^2 \right]&\leq\left\|\mathbf{x}_{\tau p}-\mathbf{x}_\star\right\|^2-2\gamma \tau
\left\langle F(\mathbf{x}_{\tau p}),\mathbf{x}_{\tau p}-\mathbf{x}_\star\right\rangle\\
&\quad+\gamma^2 \tau^2\left(\ell+\frac{2(\tau-1)}{\tau}
L_{\max}\sqrt{\frac{\kappa}{3}}\right)\left\langle F(\mathbf{x}_{\tau p}),
\mathbf{x}_{\tau p}-\mathbf{x}_\star\right\rangle\\
&\quad+\gamma^2\sigma^2\tau(\tau-1)\left(2\gamma L_{\max}(\tau+1)
+\frac{L_{\max}}{2\tau\sqrt{\ell\mu/3}}+
\frac{\gamma L_{\max}^2(\tau+1)^2}{2\tau\sqrt{\ell\mu/3}}\right)\\
&\quad+\gamma^2\sigma^2\tau\\
&\leq\left\|\mathbf{x}_{\tau p}-\mathbf{x}_\star\right\|^2-2\gamma\tau
\left\langle F(\mathbf{x}_{\tau p}),\mathbf{x}_{\tau p}-\mathbf{x}_\star
\right\rangle\\
&\quad+\gamma^2\tau^2\left(\ell+\frac{2(\tau-1)}{\tau}
L_{\max}\sqrt{\frac{\kappa}{3}}\right)\left\langle F(\mathbf{x}_{\tau p}),
\mathbf{x}_{\tau p}-\mathbf{x}_\star\right\rangle\\
&\quad+\gamma^2\sigma^2\tau(\tau-1)\left(4\gamma \tau L_{\max}
+\frac{L_{\max}}{2\tau\sqrt{\ell\mu/3}}+
\frac{2\gamma\tau L_{\max}^2}{\sqrt{\ell\mu/3}}\right)\\
&\quad+\gamma^2\sigma^2\tau,
\end{align*}
where for the last inequality, we use $\tau + 1 \leq 2 \tau$.
\end{proof}

\subsection{Proof of Lemma~\ref{lemma:methods-satisfy-key-assumption}}

Before giving the proof, we restate the lemma.

\begin{lemma}
Assume \hyperref[assumption:convexity]{\textbf{\emph{(CVX)}}},
\hyperref[assumption:smoothness]{\textbf{\emph{(SM)}}},
\hyperref[assumption:quasi-strong-monotonicity]{\textbf{\emph{(QSM)}}},
and \hyperref[assumption:star-cocoercivity]{\textbf{\emph{(SCO)}}}
hold. Let \(\kappa=\ell/\mu\), \(\sigma^2=\sum_{i=1}^n\sigma_i^2\), and
\begin{center}
\footnotesize
\(\displaystyle
\begin{aligned}
A_\tau^{\mathrm{GD}}
&:=\ell+\frac{2(\tau-1)}{\tau}L_{\max}\sqrt{\frac{\kappa}{3}},\\
B_\tau^{\mathrm{GD}} &:=\frac{\sigma^2}{\tau} \left[1+(\tau-1)L_{\max}\left(4\gamma\tau+\frac{1+4\gamma\tau^2L_{\max}}
{2\tau\sqrt{\ell\mu/3}}\right)\right],\\
U_\tau &:=\frac{\sigma^2}{\tau} \left(1+2(\tau-1)(\tau+1)\gamma L_{\max}\right),\\
A_\tau^{\mathrm{VR}} &:=4\ell+\frac{4(\tau-1)}{\tau} L_{\max}\sqrt{\frac{\kappa}{3}}.
\end{aligned}
\)
\end{center}
Then the following instances of Algorithm~\ref{alg:unified-pearl} satisfy Assumption~\ref{assm:key_assumption}.
\begin{itemize}[leftmargin=0.8em]

\item \textbf{\algname{PEARL-SGD}.}
If \hyperref[assumption:noise]{\textbf{\emph{(BV)}}} holds and
\(0<\gamma_k\equiv\gamma\leq
\frac{1}{\ell\tau+2(\tau-1)L_{\max}\sqrt{\kappa}}\),
then
\begin{align*}
    A_\tau=A_\tau^{\mathrm{GD}} \text{ and } B_\tau=B_\tau^{\mathrm{GD}}.
\end{align*}

\item \textbf{\algname{PEARL-SVRG}.}
If Assumption~\ref{assumption:finite-sum-smoothness} holds and
\(0<\gamma_k\equiv\gamma\leq\frac{1}{2L_{\max}\tau}\),
then
\begin{align*}
    A_\tau=A_\tau^{\mathrm{VR}} \text{ and } B_\tau=0.
\end{align*}

\item \textbf{\algname{C-PEARL-SGD}.}
If \hyperref[assumption:noise]{\textbf{\emph{(BV)}}} and
Assumption~\ref{assumption:unbiased-compression} hold, and
\(0<\gamma_k\equiv\gamma\leq
\frac{1}{\ell\tau+2(\tau-1)L_{\max}\sqrt{\kappa}}\),
then
\begin{align*}
    A_\tau=A_\tau^{\mathrm{GD}}+\omega\ell \text{ and } B_\tau=B_\tau^{\mathrm{GD}}+\omega U_\tau.
\end{align*}

\item \textbf{\algname{C-PEARL-SVRG}.}
If Assumptions~\ref{assumption:finite-sum-smoothness}
and~\ref{assumption:unbiased-compression} hold, and
\(0<\gamma_k\equiv\gamma\leq\frac{1}{2L_{\max}\tau}\),
then
\begin{align*}
    A_\tau=A_\tau^{\mathrm{VR}}+4\omega\ell \text{ and } B_\tau=0.
\end{align*}
\end{itemize}
\end{lemma}
\subsubsection{\algname{PEARL-SGD}}

\begin{proof}
\algname{PEARL-SGD} is obtained by choosing \(\mathcal Q^i=\mathrm{Id}_{d_i}\) and the
sample-gradient estimator $g_k^i=\nabla f_{i,\xi_k^i}(\widetilde{x}_k^i;x_{\tau p}^{-i})$ in Generalized \algname{PEARL-SGD} (\cref{alg:unified-pearl}).
Thus for every player $i \in [n]$, $q_p^i=u_p^i=x_{\tau(p+1)}^i-x_{\tau p}^i.$   Since \(\widetilde{x}_{\tau p}^i=x_{\tau p}^i\), the reconstructed synchronized strategy satisfies
\[
x_{\tau(p+1)}^i =x_{\tau p}^i+q_p^i=x_{\tau p}^i+\widetilde{x}_{\tau(p+1)}^i-\widetilde{x}_{\tau p}^i=\widetilde{x}_{\tau(p+1)}^i.
\]
Therefore, the synchronized iterates coincide with the
\algname{PEARL-SGD} iterates. From Lemma \ref{lemma:pearl-sgd-one-round-recursion}, we have an inequality of the form 
\begin{align} \notag
    \mathbb{E}_p \left[\sqnorm{\vx_{\tau(p+1)}-\vx_\star} \right] &\le\sqnorm{\vx_{\tau p}-\vx_\star}
    - 2\gamma \tau \left\langle F(\mathbf{x}_{\tau p}),\mathbf{x}_{\tau p}-\mathbf{x}_\star\right\rangle\\
    &+\gamma^2\tau^2\Bigl(A_\tau \left\langle F(\mathbf{x}_{\tau p}),\mathbf{x}_{\tau p}-\mathbf{x}_\star\right\rangle + B_\tau\Bigr).
    \label{eqn:unified-pearl-round-recursion}
\end{align}
Therefore Assumption~\ref{assm:key_assumption} holds with $A_\tau=A_\tau^{\mathrm{GD}}$ and $B_\tau =B_\tau^{\mathrm{GD}}$.
\end{proof}

\subsubsection{\algname{PEARL-SVRG}}
To verify that \algname{PEARL-SVRG} satisfies Assumption~\ref{assm:key_assumption}, we first present two lemmas that are required in our proof. For a fixed communication round \(p\),
\algname{PEARL-SVRG} uses the local estimator
\[
g_k^i=\nabla f_{i,\xi_k^i}(\widetilde{x}_k^i;x_{\tau p}^{-i})-\nabla f_{i,\xi_k^i}(x_{\tau p}^i;x_{\tau p}^{-i})+\nabla f_i(x_{\tau p}^i;x_{\tau p}^{-i}),
\]
where \(\xi_k^i\) is a mini-batch. Lemma \ref{lemma:pearl-svrg-estimator} proves that the local estimator \(g_k^i\) is unbiased conditioned on $\widetilde x_k^i$ and $ \cF_p$, and Lemma \ref{lemma:pearl-svrg-iterate-difference} bounds the local difference terms generated within one communication round. 

Note that for a fixed communication round \(p\), we use \(\widetilde{x}_k^i\) to denote the strategy of player
\(i\) during the \(\tau\) local steps. Therefore, $\widetilde{x}_{\tau p}^i=x_{\tau p}^i,$ and for \(k=\tau p,\ldots,\tau(p+1)-1\),
\[
    \widetilde{x}_{k+1}^i=\widetilde{x}_k^i-\gamma g_k^i.
\]
We use this notation only for the local strategy of player $i$ inside the round. The synchronized strategy vector is still denoted by \(\mathbf{x}_{\tau p}\).

\begin{lemma}
\label{lemma:pearl-svrg-estimator}
Under Assumption~\ref{assumption:finite-sum-smoothness}, for every
communication round \(p\), every player \(i\in[n]\), and every
\(k=\tau p,\ldots,\tau(p+1)-1\),
\begin{align}
\label{eq:svrg-unbiased-estimator}
     \mathbb E\left[g_k^i\,\middle|\,\widetilde{x}_k^i, \cF_p \right] =
    \nabla f_i\left(\widetilde{x}_k^i;x_{\tau p}^{-i}\right),
\end{align}
and
\begin{align}
\label{eq:svrg-unbiased-estimator-second-moment}
    \mathbb E\left[\left\|g_k^i\right\|^2\,\middle|\,\widetilde{x}_k^i, \cF_p \right] & \leq 2\left\|\nabla f_i\left(x_{\tau p}^i;x_{\tau p}^{-i}\right)\right\|^2 +2L_i^2\left\|\widetilde{x}_k^i-x_{\tau p}^i\right\|^2.
\end{align}
\end{lemma}

\begin{proof}

The unbiasedness in \eqref{eq:svrg-unbiased-estimator} follows from $\mathbb E\left[
\nabla f_{i,\xi_k^i}(\widetilde{x}_k^i;x_{\tau p}^{-i}) \,\middle|\,\widetilde{x}_k^i,\mathcal F_p
\right] =\nabla f_i(\widetilde{x}_k^i;x_{\tau p}^{-i})$ and $\mathbb E\left[\nabla f_{i,\xi_k^i}(x_{\tau p}^i;x_{\tau p}^{-i})
\,\middle|\,\widetilde{x}_k^i,\mathcal F_p\right]=\nabla f_i(x_{\tau p}^i;x_{\tau p}^{-i}).$
Next, for \eqref{eq:svrg-unbiased-estimator-second-moment},
\begin{align*}
    \mathbb E\left[\left\|g_k^i\right\|^2\,\middle|\,\widetilde{x}_k^i, \cF_p\right] &\leq 2\mathbb E\left[\left\|\nabla f_{i,\xi_k^i}\left(\widetilde{x}_k^i;x_{\tau p}^{-i}\right)-\nabla f_{i,\xi_k^i}\left(x_{\tau p}^i;x_{\tau p}^{-i}\right)\right\|^2\,\middle|\,\widetilde{x}_k^i,\cF_p \right] + 2\left\|\nabla f_i\left(x_{\tau p}^i;x_{\tau p}^{-i}\right)\right\|^2\\
    &\leq 2L_i^2\left\|\widetilde{x}_k^i-x_{\tau p}^i\right\|^2
    +2\left\|\nabla f_i\left(x_{\tau p}^i;x_{\tau p}^{-i}\right)\right\|^2
\end{align*}
where the last inequality uses Assumption~\ref{assumption:finite-sum-smoothness}.
\end{proof}

\begin{lemma}
\label{lemma:pearl-svrg-iterate-difference}
Under Assumption~\ref{assumption:finite-sum-smoothness}, suppose $ 0<\gamma\leq\frac{1}{2L_i\tau}.$ Then for every fixed round \(p\), player \(i\), and \(t=0,\ldots,\tau\), we have
\begin{align*}
    \mathbb E_p \left[\left\|\widetilde{x}_{\tau p+t}^i-x_{\tau p}^i\right\|^2
    \right]
    \leq4\gamma^2t^2\left\|\nabla f_i(x_{\tau p}^i;x_{\tau p}^{-i})\right\|^2.
\end{align*}
\end{lemma}

\begin{proof}
    Fix a communication round \(p\) and a player \(i\). Since $\widetilde{x}_{\tau p+t}^i-x_{\tau p}^i=-\gamma\sum_{s=0}^{t-1}g_{\tau p+s}^i,$ by using Jensen's inequality,
\begin{align}
    \mathbb E_p \left[\left\|\widetilde{x}_{\tau p+t}^i-x_{\tau p}^i\right\|^2 \right]
    & = \gamma^2\mathbb E\left[\left\|\sum_{s=0}^{t-1}g_{\tau p+s}^i\right\|^2\,\middle|\,\cF_p\right]
    \leq \gamma^2t\sum_{s=0}^{t-1}\mathbb E\left[\left\|g_{\tau p+s}^i\right\|^2\,\middle|\,\cF_p\right].
    \label{eq:pearl-svrg-local-difference-jensen}
\end{align}
By Lemma~\ref{lemma:pearl-svrg-estimator} and the tower property,
\begin{align*}
 \mathbb E\left[\left\|g_{\tau p+s}^i\right\|^2\,\middle|\, \cF_p \right] \leq
    2\left\|\nabla f_i\left(x_{\tau p}^i;x_{\tau p}^{-i}\right)\right\|^2+2L_i^2\mathbb E\left[\left\|\widetilde{x}_{\tau p+s}^i-x_{\tau p}^i\right\|^2\,\middle|\, \cF_p \right].    
\end{align*}
Substituting this into \eqref{eq:pearl-svrg-local-difference-jensen} we obtain
\begin{align}
    \mathbb E_p \left[\left\|\widetilde{x}_{\tau p+t}^i-x_{\tau p}^i\right\|^2 \right] \leq 2\gamma^2t^2\left\|\nabla f_i\left(x_{\tau p}^i;x_{\tau p}^{-i}\right)\right\|^2+2\gamma^2L_i^2t\sum_{s=0}^{t-1}\mathbb E\left[\left\|\widetilde{x}_{\tau p+s}^i-x_{\tau p}^i\right\|^2\,\middle|\,\cF_p \right].
    \label{eq:pearl-svrg-local-difference-recursion}
\end{align}
We now prove the claim by induction on \(t\). For \(t=0\), the result is
immediate because \(\widetilde{x}_{\tau p}^i=x_{\tau p}^i\). Assume that,
for all \(s=0,\ldots,t-1\),
\[
    \mathbb E_p \left[\left\|\widetilde{x}_{\tau p+s}^i-x_{\tau p}^i\right\|^2 \right] \leq 4\gamma^2s^2\left\|\nabla f_i\left(x_{\tau p}^i;x_{\tau p}^{-i}\right)\right\|^2.
\]
Then \eqref{eq:pearl-svrg-local-difference-recursion} gives
\begin{align*}
    \mathbb E_p \left[\left\|\widetilde{x}_{\tau p+t}^i-x_{\tau p}^i\right\|^2 \right] &\leq 2\gamma^2t^2\left\|\nabla f_i\left(x_{\tau p}^i;x_{\tau p}^{-i}\right)\right\|^2+8\gamma^4L_i^2t\sum_{s=0}^{t-1}s^2\left\|\nabla f_i\left(x_{\tau p}^i;x_{\tau p}^{-i}\right)\right\|^2\\
    &\leq 2\gamma^2t^2\left\|\nabla f_i\left(x_{\tau p}^i;x_{\tau p}^{-i}\right)\right\|^2+\frac{8}{3}\gamma^4L_i^2t^4\left\|\nabla f_i\left(x_{\tau p}^i;x_{\tau p}^{-i}\right)\right\|^2,
\end{align*}
where we used \(\sum_{s=0}^{t-1}s^2\leq t^3/3\). Since \(t\leq\tau\) and
\(\gamma\leq1/(2L_i\tau)\), we have $\gamma^2L_i^2t^2\leq\gamma^2L_i^2\tau^2\leq\frac14.$
Therefore, $\frac{8}{3}\gamma^4L_i^2t^4\leq\frac{2}{3}\gamma^2t^2$ and 
\begin{align*}
    \mathbb E_p \left[\left\|\widetilde{x}_{\tau p+t}^i-x_{\tau p}^i\right\|^2\right]
    & \leq
    \left(2+\frac23\right)\gamma^2t^2\left\|\nabla f_i\left(x_{\tau p}^i;x_{\tau p}^{-i}\right)\right\|^2\\
    & \leq 4\gamma^2t^2\left\|\nabla f_i\left(x_{\tau p}^i;x_{\tau p}^{-i}\right)\right\|^2.
\end{align*}
This closes the induction and proves the claim.
\end{proof}
We now verify that \algname{PEARL-SVRG} satisfies Assumption \ref{assm:key_assumption}.
\begin{proof}
For \algname{PEARL-SVRG}, since \(g_k^i\) is an unbiased estimator of
\(\nabla f_i(\widetilde{x}_k^i;x_{\tau p}^{-i})\) by Lemma \ref{lemma:pearl-svrg-estimator}, the proof of Lemma \ref{lemma:proof-outline-basic-bound} directly applies.
Hence, for any \(\alpha>0\),
\begin{align}
\label{eqn:pearl-svrg-proof-basic-bound}
    \mathbb E_p \left[\left\|\mathbf{x}_{\tau(p+1)}-\mathbf{x}_\star\right\|^2
    \right]&\leq\left(1+(\tau-1)\alpha\gamma\right)
    \left\|\mathbf{x}_{\tau p}-\mathbf{x}_\star\right\|^2-2\gamma\tau
    \left\langle F(\mathbf{x}_{\tau p}),\mathbf{x}_{\tau p}-\mathbf{x}_\star\right\rangle\\\notag
    &\quad +\frac{\gamma L_{\max}^2}{\alpha} \sum_{j=\tau p+1}^{\tau(p+1)-1}
    \mathbb E_p \left[\left\|\mathbf{x}_{\tau p}-\widetilde{\mathbf{x}}_j\right\|^2
    \right] +\mathbb E_p \left[ \left\|\mathbf{x}_{\tau p}-\widetilde{\mathbf{x}}_{\tau(p+1)}\right\|^2 \right].
\end{align}
Here
\(\widetilde{\mathbf{x}}_j=(\widetilde{x}_j^1,\ldots,\widetilde{x}_j^n)\). By \cref{lemma:pearl-svrg-iterate-difference}, for \(t=1,\ldots,\tau\),
\begin{align*}
    \mathbb E_p \left[\left\|\mathbf{x}_{\tau p}-\widetilde{\mathbf{x}}_{\tau p+t}\right\|^2 
    \right] & = \sum_{i=1}^n\mathbb E_p\left[\left\|\widetilde{x}_{\tau p+t}^i-x_{\tau p}^i\right\|^2\,\right] \\
    & \leq 4\gamma^2t^2\left\|F(\mathbf{x}_{\tau p})\right\|^2. 
\end{align*}
Substituting this into \eqref{eqn:pearl-svrg-proof-basic-bound} and using the elementary fact $\sum_{t=1}^{\tau-1}t^2 =\frac{(\tau-1)\tau(2\tau-1)}{6}\leq \frac{(\tau-1)\tau^2}{3}$ we obtain
\begin{align*}
    \mathbb E_p \left[\left\|\mathbf{x}_{\tau(p+1)}-\mathbf{x}_\star\right\|^2 \right]&\leq
    \left(1+(\tau-1)\alpha\gamma\right) \left\|\mathbf{x}_{\tau p}-\mathbf{x}_\star\right\|^2-2\gamma \tau\left\langle F(\mathbf{x}_{\tau p}),\mathbf{x}_{\tau p}-\mathbf{x}_\star\right\rangle \\
    &\quad+\left(4\gamma^2\tau^2+\frac{4\gamma^3L_{\max}^2\tau^2(\tau-1)}{3\alpha}\right)\left\|F(\mathbf{x}_{\tau p})\right\|^2 \\
    &=\left\|\mathbf{x}_{\tau p}-\mathbf{x}_\star\right\|^2-2\gamma \tau
    \left\langle F(\mathbf{x}_{\tau p}),\mathbf{x}_{\tau p}-\mathbf{x}_\star\right\rangle \\
    &\quad+(\tau-1)\alpha\gamma
    \left\|\mathbf{x}_{\tau p}-\mathbf{x}_\star\right\|^2 \\
    &\quad+\left(4\gamma^2\tau^2+\frac{4\gamma^3L_{\max}^2\tau^2(\tau-1)}{3\alpha}\right)\left\|F(\mathbf{x}_{\tau p})\right\|^2 \\
    &\hspace{-1.6em}\overset{\hyperref[assumption:quasi-strong-monotonicity]{\textbf{\emph{(QSM)}}},
    \hyperref[assumption:star-cocoercivity]{\textbf{\emph{(SCO)}}}}{\leq}
    \left\|\mathbf{x}_{\tau p}-\mathbf{x}_\star\right\|^2
    -2\gamma \tau\left\langle F(\mathbf{x}_{\tau p}),\mathbf{x}_{\tau p}-\mathbf{x}_\star\right\rangle \\
    &\quad+\frac{(\tau-1)\alpha\gamma}{\mu} \left\langle F(\mathbf{x}_{\tau p}),\mathbf{x}_{\tau p}-\mathbf{x}_\star\right\rangle \\
    &\quad+ \ell\left(4\gamma^2\tau^2+
    \frac{4\gamma^3L_{\max}^2\tau^2(\tau-1)}{3\alpha}\right)
    \left\langle F(\mathbf{x}_{\tau p}),\mathbf{x}_{\tau p}-\mathbf{x}_\star\right\rangle \\
    &=\left\|\mathbf{x}_{\tau p}-\mathbf{x}_\star\right\|^2-2\gamma \tau
    \left\langle F(\mathbf{x}_{\tau p}),\mathbf{x}_{\tau p}-\mathbf{x}_\star\right\rangle \\
    &\quad+\left(\frac{(\tau-1)\alpha\gamma}{\mu}+4\ell\gamma^2\tau^2
    +\frac{4\ell\gamma^3L_{\max}^2\tau^2(\tau-1)}{3\alpha}\right)\\
    &\hspace{15em}\left\langle F(\mathbf{x}_{\tau p}),\mathbf{x}_{\tau p}-\mathbf{x}_\star\right\rangle\\
    &\leq \left\|\mathbf{x}_{\tau p}-\mathbf{x}_\star\right\|^2
    -2\gamma \tau\left\langle F(\mathbf{x}_{\tau p}),\mathbf{x}_{\tau p}-\mathbf{x}_\star\right\rangle\\
    &\quad+\gamma^2 \tau^2A_\tau(\alpha)\left\langle F(\mathbf{x}_{\tau p}),\mathbf{x}_{\tau p}-\mathbf{x}_\star\right\rangle,
\end{align*}
where $A_\tau(\alpha):=\frac{(\tau-1)\alpha}{\mu\gamma\tau^2}+4\ell+ \frac{4\ell\gamma L_{\max}^2(\tau-1)}{3\alpha}.$
Now we choose $\alpha= 2\gamma\tau L_{\max}\sqrt{\frac{\ell\mu}{3}}$ which minimizes $A_\tau(\alpha)$. With this choice of $\alpha$, the previous inequality becomes
\begin{align}
\label{eqn:pearl-svrg-descent-bound}
    \mathbb E_p \left[\left\|\mathbf{x}_{\tau(p+1)}-\mathbf{x}_\star\right\|^2
    \right] \leq\left\|\mathbf{x}_{\tau p}-\mathbf{x}_\star\right\|^2
    -2\gamma \tau\left\langle F(\mathbf{x}_{\tau p}),\mathbf{x}_{\tau p}-\mathbf{x}_\star\right\rangle +\gamma^2 \tau^2A_\tau^{\mathrm{VR}}
    \left\langle F(\mathbf{x}_{\tau p}),\mathbf{x}_{\tau p}-\mathbf{x}_\star\right\rangle.
\end{align}
where \[
    A_\tau(\alpha)=4\ell+\frac{4(\tau-1)}{\tau}L_{\max}\sqrt{\frac{\kappa}{3}} :=
    A_\tau^{\mathrm{VR}} . 
\]
Thus Assumption~\ref{assm:key_assumption} holds with \(A_\tau=A_\tau^{\mathrm{VR}}\), \(B_\tau=0\).
\end{proof}

\subsubsection{\algname{C-PEARL-SGD}}
\begin{proof}
For each player \(i\in[n]\), let \(\widetilde{x}^i_{\tau(p+1)}\) denote the strategy of player \(i\) after the \(\tau\) local steps in communication round \(p\), before compression. Thus, by the algorithm, $u_p^i=\widetilde{x}^i_{\tau(p+1)} - \widetilde{x}_{\tau p}^i=\widetilde{x}^i_{\tau(p+1)} - x_{\tau p}^i,$ and $q_p^i = \mathcal Q^i(u_p^i).$ The master concatenates and broadcasts $\mathbf q_p=(q_p^1,\ldots,q_p^n),$ and each player reconstructs the synchronized strategy vector $\mathbf{x}_{\tau(p+1)}= \mathbf{x}_{\tau p}+\mathbf q_p.$ 

Let $\widetilde{\mathbf{x}}_{\tau(p+1)}=(\widetilde{x}_{\tau(p+1)}^1,\ldots,
\widetilde{x}_{\tau(p+1)}^n).$ Since $\widetilde{\mathbf{x}}_{\tau(p+1)}=\mathbf{x}_{\tau p}+\mathbf u_p,$ we have $\mathbf{x}_{\tau(p+1)}-\mathbf{x}_\star=\widetilde{\mathbf{x}}_{\tau(p+1)}-\mathbf{x}_\star
+\mathbf q_p-\mathbf u_p.$
Therefore,
\begin{align}
    \left\|\mathbf{x}_{\tau(p+1)}-\mathbf{x}_\star\right\|^2=\left\|\widetilde{\mathbf{x}}_{\tau(p+1)}
    -\mathbf{x}_\star\right\|^2+2\left\langle\widetilde{\mathbf{x}}_{\tau(p+1)}-\mathbf{x}_\star,
    \mathbf q_p-\mathbf u_p\right\rangle + \left\|\mathbf q_p-\mathbf u_p\right\|^2.
    \label{eq:cpearl-compression-expansion}
\end{align}
By Assumption~\ref{assumption:unbiased-compression}, $\mathbb E_\cQ \left[ \mathbf q_p-\mathbf u_p \right] = \mathbb E \left[\mathbf q_p-\mathbf u_p\,\middle|\,\mathcal F_p,\mathbf u_p\right] = 0$ and
\[ \mathbb E_\cQ \left[\left\|\mathbf q_p-\mathbf u_p\right\|^2 \right]
= \mathbb E\left[\left\|\mathbf q_p-\mathbf u_p\right\|^2
\,\middle|\,\mathcal F_p,\mathbf u_p
\right]\leq\omega\left\|\mathbf u_p\right\|^2.
\]
Moreover,
\(\widetilde{\mathbf{x}}_{\tau(p+1)}-\mathbf{x}_\star
=\mathbf{x}_{\tau p}+\mathbf u_p-\mathbf{x}_\star\)
is deterministic conditioned on $\cF_p$ and \(\mathbf u_p\); hence, taking conditional expectation in \eqref{eq:cpearl-compression-expansion} given $\cF_p$ and $\mathbf u_p$, and applying the tower property, we obtain
\begin{align}
\label{eq:cpearl-after-compression}
    & \mathbb E_p\left[ \left\|\mathbf{x}_{\tau(p+1)}-\mathbf{x}_\star\right\|^2 \right] = 
    \mathbb E_p \left[ \mathbb E_\cQ \left[ \left\|\mathbf{x}_{\tau(p+1)}-\mathbf{x}_\star\right\|^2 \right] \right] \notag \\
    & \leq  \mathbb E_p\left[\left\|\widetilde{\mathbf{x}}_{\tau(p+1)}-\mathbf{x}_\star\right\|^2
    \right] + \omega\mathbb E_p\left[\left\|\mathbf u_p\right\|^2 \right].
\end{align}

We first bound the first term in
\eqref{eq:cpearl-after-compression}. Since
\(\widetilde{\mathbf{x}}_{\tau(p+1)}\) is the strategy after the
\(\tau\) local \algname{PEARL-SGD} steps before compression,
Lemma~\ref{lemma:pearl-sgd-one-round-recursion} gives
\begin{align}
\label{eq:cpearl-local-strategy-after-tau-steps}
    \mathbb E_p\left[\left\|\widetilde{\mathbf{x}}_{\tau(p+1)}-\mathbf{x}_\star\right\|^2
    \right] &\leq\left\|\mathbf{x}_{\tau p}-\mathbf{x}_\star\right\|^2-2\gamma\tau
    \left\langle F(\mathbf{x}_{\tau p}),\mathbf{x}_{\tau p}-\mathbf{x}_\star \right\rangle \notag\\
    &\quad +\gamma^2\tau^2 \left(A_\tau^{\mathrm{GD}}\left\langle
    F(\mathbf{x}_{\tau p}),\mathbf{x}_{\tau p}-\mathbf{x}_\star\right\rangle
    +B_\tau^{\mathrm{GD}} \right).
\end{align}

Next, using Lemma~\ref{lemma:local-GDA-iterate-difference-bound-stochastic-case}
with \(t=\tau\), summing over \(i=1,\ldots,n\), and using
\(L_i\leq L_{\max}\), we get
\begin{align}
    \mathbb E_p\left[\left\|\mathbf u_p\right\|^2\right]&=\mathbb E_p\left[\sum_{i=1}^n\left\|\widetilde{x}_{\tau(p+1)}^i-x_{\tau p}^i
    \right\|^2\right] \notag\\
    &\leq\gamma^2\tau^2\sum_{i=1}^n\left\|\nabla f_i(x_{\tau p}^i;x_{\tau p}^{-i})
    \right\|^2 + \gamma^2\tau\sum_{i=1}^n\left(1+2(\tau-1)(\tau+1)\gamma L_i\right)\sigma_i^2 \notag\\
    &\leq\gamma^2\tau^2 \left\|F(\mathbf{x}_{\tau p})\right\|^2+\gamma^2\tau\left(1+2(\tau-1)(\tau+1)\gamma L_{\max}\right)\sigma^2 \notag\\
    &=\gamma^2\tau^2\left( \left\|F(\mathbf{x}_{\tau p})\right\|^2+ U_\tau\right) \notag\\
    &\overset{\hyperref[assumption:star-cocoercivity]{\textbf{\emph{(SCO)}}}}{\leq}
    \gamma^2\tau^2\left(\ell\left\langle F(\mathbf{x}_{\tau p}),
    \mathbf{x}_{\tau p}-\mathbf{x}_\star\right\rangle
    +U_\tau\right),
    \label{eq:cpearl-playerwise-update-bound}
\end{align}
where $U_\tau := \frac{\sigma^2}{\tau}\left(1+2(\tau-1)(\tau+1)\gamma L_{\max}\right).$
Substituting \eqref{eq:cpearl-local-strategy-after-tau-steps} and
\eqref{eq:cpearl-playerwise-update-bound} into \eqref{eq:cpearl-after-compression} yields,
\begin{align*}
    \mathbb E_p\left[\left\|\mathbf{x}_{\tau(p+1)}-\mathbf{x}_\star\right\|^2\right]&\leq \left\|\mathbf{x}_{\tau p}-\mathbf{x}_\star\right\|^2-2\gamma\tau\left\langle F(\mathbf{x}_{\tau p}),\mathbf{x}_{\tau p}-\mathbf{x}_\star\right\rangle\\
    &\quad+\gamma^2\tau^2 \left(
    A_\tau^{\mathrm{GD}}
    \left\langle F(\mathbf{x}_{\tau p}),\mathbf{x}_{\tau p}-\mathbf{x}_\star\right\rangle
    +B_\tau^{\mathrm{GD}}
    \right)\\
    &\quad+\omega\gamma^2\tau^2 \left(
    \ell \left\langle F(\mathbf{x}_{\tau p}),\mathbf{x}_{\tau p}-\mathbf{x}_\star\right\rangle
    +U_\tau
    \right)\\
    &= \left\|\mathbf{x}_{\tau p}-\mathbf{x}_\star\right\|^2
    -2\gamma\tau\left\langle F(\mathbf{x}_{\tau p}),\mathbf{x}_{\tau p}-\mathbf{x}_\star\right\rangle\\
    &\quad+\gamma^2\tau^2 \left(
    \left(A_\tau^{\mathrm{GD}}+\omega\ell\right)
    \left\langle F(\mathbf{x}_{\tau p}),\mathbf{x}_{\tau p}-\mathbf{x}_\star\right\rangle
    +B_\tau^{\mathrm{GD}}
    +\omega U_\tau
    \right).
\end{align*}
This is exactly \eqref{eq:unified_round_descent} with $A_\tau=A_\tau^{\mathrm{GD}}+\omega\ell$, and $ B_{\tau}=B_\tau^{\mathrm{GD}}+\omega U_\tau.$ Therefore, \algname{C-PEARL-SGD} satisfies Assumption~\ref{assm:key_assumption}.
\end{proof}
\subsubsection{\algname{C-PEARL-SVRG}}
\begin{proof}
For \(k=\tau p,\ldots,\tau(p+1)-1\), the local update before compression is
the same as in \algname{PEARL-SVRG}. For each player \(i\in[n]\), $\widetilde{x}_{\tau p}^i=x_{\tau p}^i,$ and $\widetilde{x}_{k+1}^i=\widetilde{x}_{k}^i-\gamma g_k^i,$
where $g_k^i=\nabla f_{i,\xi_k^i}\left(\widetilde{x}_k^i;x_{\tau p}^{-i}\right)
-\nabla f_{i,\xi_k^i}\left(x_{\tau p}^i;x_{\tau p}^{-i}\right)
+\nabla f_i\left(x_{\tau p}^i;x_{\tau p}^{-i}\right).$
After the \(\tau\) local steps, player \(i\) forms $u_p^i
=\widetilde{x}_{\tau(p+1)}^i-x_{\tau p}^i,$ and compresses it as $q_p^i=\mathcal Q^i(u_p^i).$
Let $\mathbf u_p=(u_p^1,\ldots,u_p^n),$ and $\mathbf q_p=(q_p^1,\ldots,q_p^n),$
and  $\widetilde{\mathbf{x}}_{\tau(p+1)}=(\widetilde{x}_{\tau(p+1)}^1,\ldots,
\widetilde{x}_{\tau(p+1)}^n)$.
Again, we have
$\mathbf{x}_{\tau(p+1)} - \mathbf{x}_\star = (\widetilde{\mathbf{x}}_{\tau(p+1)} - \mathbf{x}_\star)
+ (\mathbf q_p-\mathbf u_p)$, so we have the same identities~\eqref{eq:cpearl-compression-expansion} and \eqref{eq:cpearl-after-compression} as in \algname{C-PEARL-SGD}.

Here, our $\widetilde{\vx}_{\tau (p+1)}$ before compression is constructed in the same way as \algname{PEARL-SVRG} conditioned on $\cF_p$, so following the same derivation for \eqref{eqn:pearl-svrg-descent-bound} we have 
\begin{align}
    \mathbb E_p\left[\left\|\widetilde{\mathbf{x}}_{\tau(p+1)}
    -\mathbf{x}_\star\right\|^2 \right] &\leq\left\|\mathbf{x}_{\tau p}-\mathbf{x}_\star\right\|^2-2\gamma\tau
    \left\langle F(\mathbf{x}_{\tau p}),\mathbf{x}_{\tau p}-\mathbf{x}_\star
    \right\rangle+\gamma^2\tau^2 A_\tau^{\mathrm{VR}}
    \left\langle F(\mathbf{x}_{\tau p}),
    \mathbf{x}_{\tau p}-\mathbf{x}_\star\right\rangle .
    \label{eq:cpearl-vr-local-descent}
\end{align}
Next, by Lemma~\ref{lemma:pearl-svrg-iterate-difference} with
\(t=\tau\), the assumed step-size, and \(L_i\le L_{\max}\),
\begin{align}
    \mathbb E_p\left[
    \left\|\mathbf u_p\right\|^2 \right]&=
    \mathbb E_p\left[
    \sum_{i=1}^n\left\|\widetilde{x}_{\tau(p+1)}^i-x_{\tau p}^i
    \right\|^2\right]\notag\\
    &\leq4\gamma^2\tau^2\sum_{i=1}^n \left\|
    \nabla f_i\left(x_{\tau p}^i;x_{\tau p}^{-i}\right)\right\|^2
    \notag\\
    &= 4\gamma^2\tau^2 \left\|F(\mathbf{x}_{\tau p})\right\|^2\notag\\
    &\overset{\hyperref[assumption:star-cocoercivity]{\textbf{\emph{(SCO)}}}}{\leq}
    4\gamma^2\tau^2\ell\left\langle
    F(\mathbf{x}_{\tau p}), \mathbf{x}_{\tau p}-\mathbf{x}_\star
    \right\rangle .
    \label{eq:cpearl-vr-update-bound}
\end{align}
Substituting \eqref{eq:cpearl-vr-local-descent} and
\eqref{eq:cpearl-vr-update-bound} into
\eqref{eq:cpearl-after-compression}, we get
\begin{align*}
    \mathbb E_p\left[\left\|\mathbf{x}_{\tau(p+1)}-\mathbf{x}_\star\right\|^2
    \right] &\leq\left\|\mathbf{x}_{\tau p}-\mathbf{x}_\star\right\|^2
    -2\gamma\tau\left\langle F(\mathbf{x}_{\tau p}),\mathbf{x}_{\tau p}-\mathbf{x}_\star \right\rangle\\
    &\quad+\gamma^2\tau^2 A_\tau^{\mathrm{VR}}\left\langle
    F(\mathbf{x}_{\tau p}),\mathbf{x}_{\tau p}-\mathbf{x}_\star
    \right\rangle\\
    &\quad+4\omega\gamma^2\tau^2\ell
    \left\langle F(\mathbf{x}_{\tau p}),
    \mathbf{x}_{\tau p}-\mathbf{x}_\star
    \right\rangle\\
    &=\left\|\mathbf{x}_{\tau p}-\mathbf{x}_\star\right\|^2-2\gamma\tau
    \left\langle F(\mathbf{x}_{\tau p}),
    \mathbf{x}_{\tau p}-\mathbf{x}_\star
    \right\rangle\\
    &\quad+\gamma^2\tau^2\left(A_\tau^{\mathrm{VR}}+4\omega\ell\right)
    \left\langle F(\mathbf{x}_{\tau p}),\mathbf{x}_{\tau p}-\mathbf{x}_\star
    \right\rangle .
\end{align*}
This is \eqref{eq:unified_round_descent} with $A_\tau=A_\tau^{\mathrm{VR}}+4\omega\ell,$ and $B_\tau=0.$ Therefore, \algname{C-PEARL-SVRG} satisfies
Assumption~\ref{assm:key_assumption}.
\end{proof}
\subsection{Proof of Theorem~\ref{thm:main_theorem}}
\begin{proof}
Starting from Assumption~\ref{assm:key_assumption} and using \hyperref[assumption:quasi-strong-monotonicity]{\textbf{\emph{(QSM)}}},
\begin{align}\notag
    \mathbb E_p\left[\left\|\mathbf{x}_{\tau(p+1)}-\mathbf{x}_\star\right\|^2\right]
    &\overset{\eqref{eq:unified_round_descent}}{\leq}\left\|\mathbf{x}_{\tau p}-\mathbf{x}_\star\right\|^2-2\gamma\tau\left\langle F(\mathbf{x}_{\tau p}),\mathbf{x}_{\tau p}-\mathbf{x}_\star\right\rangle\\\notag
    &\quad+\gamma^2\tau^2\left(A_\tau\left\langle F(\mathbf{x}_{\tau p}),\mathbf{x}_{\tau p}-\mathbf{x}_\star\right\rangle+B_{\tau}\right)\\\notag
    &= \left\|\mathbf{x}_{\tau p}-\mathbf{x}_\star\right\|^2-\gamma\tau\bigl(2-\gamma\tau A_\tau\bigr)\left\langle F(\mathbf{x}_{\tau p}),\mathbf{x}_{\tau p}-\mathbf{x}_\star\right\rangle+\gamma^2\tau^2B_{\tau}\\\notag
    &\hspace{-1.0em}\overset{\hyperref[assumption:quasi-strong-monotonicity]{\textbf{\emph{(QSM)}}}}{\leq} \left\|\mathbf{x}_{\tau p}-\mathbf{x}_\star\right\|^2-\gamma\tau\mu\bigl(2-\gamma\tau A_\tau\bigr)
    \left\|\mathbf{x}_{\tau p}-\mathbf{x}_\star\right\|^2+\gamma^2\tau^2B_{\tau}\notag\\
    &= \left(1-\gamma\tau\mu\zeta\right)
    \left\|\mathbf{x}_{\tau p}-\mathbf{x}_\star\right\|^2
    +\gamma^2\tau^2B_\tau.
    \label{eq:main_proof_1}
\end{align}
Since $A_\tau\geq\mu$ and $0<\gamma\tau A_\tau<2$, we have $ 0<
\gamma\tau\mu\zeta=\frac{\mu}{A_\tau}\left(\gamma\tau A_\tau\right)\left(2-\gamma\tau A_\tau\right)\leq 1.$ Taking total expectation and using the tower property, we obtain
\begin{align*}
    \mathbb E[\left\|\mathbf{x}_{\tau (p+1)}-\mathbf{x}_\star\right\|^2]&\leq(1-\gamma\tau\mu\zeta)\mathbb E[\left\|\mathbf{x}_{\tau p}-\mathbf{x}_\star\right\|^2] + \gamma^2\tau^2B_{\tau}.
\end{align*}
Finally, recursive application of the above inequality for \(p=0,\ldots,R-1\) yields
\begin{align*}
    \mathbb E[\left\|\mathbf{x}_{\tau R}-\mathbf{x}_\star\right\|^2]
    &\leq(1-\gamma\tau\mu\zeta)^R \|\mathbf{x}_{0}-\mathbf{x}_\star\|^2 +\frac{\gamma\tau B_{\tau}}{\mu\zeta}.
\end{align*}
\end{proof}

\section{Comparison with Algorithms in Distributed Optimization}
\label{appendix:details_exact_convergence_increment_sharing}
In Section~\ref{section:MpFL-concept}, we compared the increment-sharing mechanism with two compression-based methods in distributed minimization: QSGD (Quantized SGD)~\cite{alistarh2017qsgd,gorbunov2020unified} and DIANA~\cite{mishchenko2025distributed}. Here, we explain this connection in more detail. Since distributed minimization and MpFL are different settings, we only compare how the methods choose the quantity to compress and how this affects the error due to compression.

\subsection{Problem Setting}
Consider the standard distributed minimization problem
\begin{align*}
    \min_{x\in\mathbb R^d}f(x):=\frac{1}{n}\sum_{i=1}^n f_i(x),
\end{align*}
where worker \(i\) stores the local objective \(f_i\). All workers optimize the same shared variable \(x\), and the server averages the information received from them. At a solution \(x_\star\), the global gradient satisfies
\(\nabla f(x_\star)=\frac{1}{n}\sum_{i=1}^n\nabla f_i(x_\star)=0\), but the individual local gradients \(\nabla f_i(x_\star)\) need not be zero.

\subsection{Comparison with QSGD and DIANA}
Figure~\ref{fig:c-pearl-qsgd-diana-comparison} compares QSGD~\cite{alistarh2017qsgd, mishchenko2025distributed}, DIANA~\cite{mishchenko2025distributed}, and \algname{C-PEARL-SGD} with the algorithms grouped by their problem setting.  For QSGD, we follow the distributed formulation in \cite[Algorithm~2]{mishchenko2025distributed}, and for DIANA, we use \cite[Algorithm~1]{mishchenko2025distributed} without momentum or regularization. In both cases, $\mathcal Q^i$ denotes the corresponding quantization map. For this comparison, we consider deterministic gradients in all three algorithms. For QSGD and DIANA, $g_p^i=\nabla f_i(x_p)$ denotes the local gradient. The server averages the compressed information $q_p^i$ received from the workers to obtain
$\overline q_p=\frac{1}{n}\sum_{i=1}^n q_p^i$. In contrast, in MpFL, each player $i$ minimizes its own objective $f_i(x^i;x^{-i})$ with respect to its strategy $x^i$. For \algname{C-PEARL-SGD}, $g_k^i=\nabla f_i(\widetilde x_k^i;x_{\tau p}^{-i})$ is the local gradient in Algorithm~\ref{alg:unified-pearl},
and the server concatenates the compressed local strategy increments $q_p^i$ to form $\mathbf q_p=(q_p^1,\ldots,q_p^n)$ rather than averaging them.

\begin{figure}[H]
\centering

\setlength{\fboxsep}{4pt}
\setlength{\fboxrule}{0.4pt}
\fbox{%
\begin{minipage}[t]{0.625\textwidth}
\vspace{0pt}
\centering

\small\textbf{Distributed Minimization}\par
\vspace{0.25em}

\noindent
\begin{minipage}[t]{0.49\linewidth}
\vspace{0pt}\raggedright
\hrule\vspace{0.35em}
\textbf{QSGD}~\cite{alistarh2017qsgd,mishchenko2025distributed}\par
\vspace{0.25em}\hrule\vspace{0.35em}
\footnotesize
\algrenewcommand\algorithmicindent{0.7em}
\algrenewcommand\alglinenumber[1]{\scriptsize #1:}
\begin{algorithmic}[1]
\State \textbf{Input:} $\gamma,R,\{\mathcal Q^i\},x_0$
\For{$p=0,\ldots,R-1$}
    \State Server broadcasts $x_p$
    \For{$i=1,\ldots,n$ in parallel}
        \State Compute $g_p^i$
        \State $q_p^i\gets\mathcal Q^i(\textcolor{red}{g_p^i})$
        \State Send $q_p^i$ to server
    \EndFor
    \State $\overline q_p\gets\frac{1}{n}\sum_{i=1}^n q_p^i$
    \State $x_{p+1}\gets x_p-\gamma\overline q_p$
\EndFor
\State \textbf{Output:} $x_R$
\end{algorithmic}
\vspace{0.35em}\hrule
\end{minipage}%
\hfill
\begin{minipage}[t]{0.49\linewidth}
\vspace{0pt}\raggedright
\hrule\vspace{0.35em}
\textbf{DIANA}~\cite{mishchenko2025distributed}\par
\vspace{0.25em}\hrule\vspace{0.35em}
\footnotesize
\algrenewcommand\algorithmicindent{0.7em}
\algrenewcommand\alglinenumber[1]{\scriptsize #1:}
\begin{algorithmic}[1]
\State \textbf{Input:} $\gamma,\alpha,R,\{\mathcal Q^i\},x_0,\{h_0^i\}$
\State $h_0\gets\frac{1}{n}\sum_{i=1}^n h_0^i$
\For{$p=0,\ldots,R-1$}
    \State Server broadcasts $x_p$
    \For{$i=1,\ldots,n$ in parallel}
        \State Compute $g_p^i$
        \State $\Delta_p^i\gets g_p^i-h_p^i$
        \State $q_p^i\gets\mathcal Q^i(\textcolor{red}{\Delta_p^i})$
        \State $h_{p+1}^i\gets h_p^i+\alpha q_p^i$
        \State Send $q_p^i$ to server
    \EndFor
    \State $\overline q_p\gets\frac{1}{n}\sum_{i=1}^n q_p^i$
    \State $g_p\gets h_p+\overline q_p$
    \State $x_{p+1}\gets x_p-\gamma g_p$
    \State $h_{p+1}\gets h_p+\alpha\overline q_p$
\EndFor
\State \textbf{Output:} $x_R$
\end{algorithmic}
\vspace{0.35em}\hrule
\end{minipage}

\end{minipage}%
}%
\hfill
\fbox{%
\begin{minipage}[t]{0.315\textwidth}
\vspace{0pt}
\centering

\small\textbf{MpFL}\par
\vspace{0.25em}

\noindent
\begin{minipage}[t]{\linewidth}
\vspace{0pt}\raggedright
\hrule\vspace{0.35em}
\textbf{\algname{C-PEARL-SGD}}\par
\vspace{0.25em}\hrule\vspace{0.35em}
\footnotesize
\algrenewcommand\algorithmicindent{0.7em}
\algrenewcommand\alglinenumber[1]{\scriptsize #1:}
\begin{algorithmic}[1]
\State \textbf{Input:} $\gamma_k,\tau,R,\{\mathcal Q^i\},\{x_0^i\}$
\State Server forms and broadcasts $\mathbf{x}_0$
\For{$p=0,\ldots,R-1$}
    \For{$i=1,\ldots,n$ in parallel}
        \State $\widetilde x_{\tau p}^i\gets x_{\tau p}^i$
        \For{$k=\tau p,\ldots,\tau(p+1)-1$}
            \State Compute $g_k^i$
            \State $\widetilde x_{k+1}^i
            \gets\widetilde x_k^i-\gamma_k g_k^i$
        \EndFor
        \State $u_p^i\gets
        \widetilde x_{\tau(p+1)}^i-x_{\tau p}^i$
        \State $q_p^i\gets
        \mathcal Q^i(\textcolor{red}{u_p^i})$
        \State Send $q_p^i$ to server
    \EndFor
    \State $\mathbf q_p\gets(q_p^1,\ldots,q_p^n)$
    \State Server broadcasts $\mathbf q_p$
    \State $\mathbf{x}_{\tau(p+1)}
    \gets\mathbf{x}_{\tau p}+\mathbf q_p$
\EndFor
\State \textbf{Output:} $\mathbf{x}_{\tau R}$
\end{algorithmic}
\vspace{0.35em}\hrule
\end{minipage}

\end{minipage}
}

\caption{Comparison of the quantities compressed by QSGD and DIANA in distributed minimization (left) and by \algname{C-PEARL-SGD} in MpFL (right). The input to each compressor is highlighted in red.}
\label{fig:c-pearl-qsgd-diana-comparison}
\end{figure}

Having described the different problem settings and update rules in Fig.~\ref{fig:c-pearl-qsgd-diana-comparison}, we now focus on the quantities that these methods compress and how this choice affects the compression error. In distributed minimization, QSGD~\cite{alistarh2017qsgd,mishchenko2025distributed} directly compresses the local gradient $g_p^i=\nabla f_i(x_p)$, which need not approach zero as the algorithm approaches the solution $x_\star$. This is because even though the global gradient vanishes at the solution, the individual local gradients $\nabla f_i(x_\star)$ need not be zero, and as a result, the error due to compressing these gradients may not vanish, which can prevent exact convergence. DIANA~\cite{mishchenko2025distributed} addresses this by compressing the difference between the local gradient and a memory vector, rather than the gradient itself. Each worker maintains a memory vector $h_p^i$ and compresses
$\Delta_p^i=g_p^i-h_p^i$. Under the convergence conditions of DIANA, both $g_p^i$ and $h_p^i$ approach $\nabla f_i(x_\star)$, so the difference $\Delta_p^i$ approaches zero~\cite{mishchenko2025distributed} and therefore the compression error also vanishes as the method approaches the solution.

The same idea explains increment sharing in MpFL. \algname{C-PEARL-SGD} compresses the local strategy increment \(u_p^i=\widetilde x_{\tau(p+1)}^i-x_{\tau p}^i\), which approaches zero as the synchronized and locally updated strategies approach the same equilibrium. Unlike DIANA, we do not need to maintain an
additional memory vector since the synchronized strategy $x_{\tau p}^i$ is already available at the beginning of each round, and the local updates naturally produce the increment $u_p^i$. Therefore, despite the differences in the setting, both DIANA and \algname{C-PEARL-SGD} use the same idea of compressing
a quantity that approaches zero, allowing the compression error to vanish as the method converges.
\end{document}